\documentclass[10pt]{article}

\usepackage{amscd,amsmath, amssymb, fancyhdr, mathbbol, dutchcal, url}
\usepackage[matrix,arrow,curve]{xy}
\usepackage[backref=page]{hyperref}

\renewcommand*{\backrefalt}[4]{%
	\ifcase #1 (Not cited.)%
	\or        (Cited on page~#2.)%
	\else      (Cited on pages~#2.)%
	\fi}

\hypersetup{
	colorlinks   = true,
	citecolor    = magenta,
	linkcolor    =blue,
	urlcolor     =magenta	
}

\numberwithin{equation}{section}

\newcommand{\version}{version 1.0,\ \ March 9, 2021}

\def\eqref#1{(\ref{#1})}

\newcommand{\g}{{\mathfrak g}}
\newcommand{\h}{{\mathfrak h}}
\newcommand{\arrow}{{\:\longrightarrow\:}}
\newcommand{\Z}{{\Bbb Z}}
\def\C{{\Bbb C}}
\newcommand{\R}{{\Bbb R}}
\newcommand{\Q}{{\Bbb Q}}
\renewcommand{\H}{{\Bbb H}}

\def\1{\sqrt{-1}\:}
\newcommand{\restrict}[1]{{\left|_{{\phantom{|}\!\!}_{#1}}\right.}}
\newcommand{\cntrct}                
{\hspace{2pt}\raisebox{1pt}{\text{$\lrcorner$}}\hspace{2pt}}
\newcommand{\calO}{\mathcal O}

\renewcommand{\t}{\mathfrak t}
\newcommand{\z}{\mathfrak z}

\renewcommand{\tilde}{\widetilde}
\renewcommand{\bar}{\overline}
\renewcommand{\phi}{\varphi}
\renewcommand{\epsilon}{\varepsilon}

\newcommand{\End}{\operatorname{End}}

\newcommand{\Id}{\operatorname{Id}}

\newcommand{\inv}{\operatorname{\it inv}}
\newcommand{\Vol}{\operatorname{Vol}}

\newcommand{\Lie}{\operatorname{Lie}}

\newcommand{\Alb}{\operatorname{Alb}}

\renewcommand{\Re}{\operatorname{Re}}

\newcounter{Mycounter}[section]
\newcounter{lemma}[section]
\renewcommand{\thelemma}{{Lemma \thesection.\arabic{lemma}}}
\newcommand{\lemma}{%
    \setcounter{lemma}{\value{Mycounter}}
    \refstepcounter{lemma}
    \stepcounter{Mycounter}
    {\noindent \bf \thelemma:\ }}

\newcounter{claim}[section]
\newcounter{sublemma}[section]
\newcounter{corollary}[section]
\renewcommand{\thecorollary}{{Corollary \thesection.\arabic{corollary}}}
\newcommand{\corollary}{%
    \setcounter{corollary}{\value{Mycounter}}
    \refstepcounter{corollary}
    \stepcounter{Mycounter}
    {\noindent \bf \thecorollary:\ }}

\newcounter{theorem}[section]
\renewcommand{\thetheorem}{{Theorem \thesection.\arabic{theorem}}}
\newcommand{\theorem}{%
    \setcounter{theorem}{\value{Mycounter}}
    \refstepcounter{theorem}
    \stepcounter{Mycounter}
    {\noindent \bf \thetheorem:\ }}

\newcounter{conjecture}[section]
\renewcommand{\theconjecture}{{Conjecture \thesection.\arabic{conjecture}}}
\newcommand{\conjecture}{%
    \setcounter{conjecture}{\value{Mycounter}}
    \refstepcounter{conjecture}
    \stepcounter{Mycounter}
    {\noindent \bf \theconjecture:\ }}

\newcounter{proposition}[section]
\renewcommand{\theproposition}
      {{Proposition \thesection.\arabic{proposition}}}
\newcommand{\proposition}{%
    \setcounter{proposition}{\value{Mycounter}}
    \refstepcounter{proposition}
    \stepcounter{Mycounter}
    {\noindent \bf \theproposition:\ }}

\newcounter{definition}[section]
\renewcommand{\thedefinition}
      {{Definition~\thesection.\arabic{definition}}}
\newcommand{\definition}{%
    \setcounter{definition}{\value{Mycounter}}
    \refstepcounter{definition}
    \stepcounter{Mycounter}
    {\noindent \bf \thedefinition:\ }}

\newcounter{example}[section]
\renewcommand{\theexample}{{Example \thesection.\arabic{example}}}
\newcommand{\example}{%
    \setcounter{example}{\value{Mycounter}}
    \refstepcounter{example}
    \stepcounter{Mycounter}
    {\noindent \bf \theexample:\ }}

\newcounter{remark}[section]
\renewcommand{\theremark}{{Remark \thesection.\arabic{remark}}}
\newcommand{\remark}{%
    \setcounter{remark}{\value{Mycounter}}
    \refstepcounter{remark}
    \stepcounter{Mycounter}
    {\noindent \bf \theremark:\ }}

\newcounter{problem}[section]
\newcounter{question}[section]
\renewcommand{\thequestion}{{Question \thesection.\arabic{question}}}
\newcommand{\question}{%
    \setcounter{question}{\value{Mycounter}}
    \refstepcounter{question}
    \stepcounter{Mycounter}
    {\noindent \bf \thequestion:\ }}

\newcommand{\proof}{\noindent{\bf Proof:\ }}

\makeatletter

\@addtoreset{equation}{section} \@addtoreset{footnote}{section}

\def\x@arrow{\DOTSB\Relbar}
\def\xlongrightarrowfill@{\arrowfill@\relbar\relbar\longrightarrow}
\newcommand{\xlongrightarrow}[2][]{%
        \ext@arrow 0099\xlongrightarrowfill@{#1}{#2}}
\makeatother

\def\blacksquare{\hbox{\vrule width 5pt height 5pt depth 0pt}}
\def\endproof{\blacksquare}

\begin{document}

\begin{center}
{\LARGE\bf
Algebraic dimension and complex subvarieties of hypercomplex nilmanifolds\\[4mm]
}

Anna Abasheva, Misha Verbitsky,\footnote{Both authors acknowledge support of 
HSE University basic research program,
A.\:A. is supported by Simons Foundation, 
M.\:V. is supported by FAPERJ E-26/202.912/2018 
and CNPq - Process 313608/2017-2. \\

}

\end{center}

{\small \hspace{0.10\linewidth}
\begin{minipage}[t]{0.85\linewidth}
{\bf Abstract.} 
A nilmanifold is a (left) quotient of a nilpotent Lie group by a cocompact lattice.
A hypercomplex structure on a manifold is a triple of complex structure operators
satisfying the quaternionic relations. A hypercomplex nilmanifold is a compact 
quotient of a nilpotent Lie group equipped with a left-invariant hypercomplex structure.
Such a manifold admits a whole 2-dimensional sphere $S^2$ of complex structures
induced by quaternions. We prove that for any hypercomplex nilmanifold $M$
and a generic complex structure $L\in S^2$, the complex
manifold $(M,L)$ has algebraic dimension 0.
A stronger result is proven when the hypercomplex nilmanifold is
abelian. Consider the Lie algebra of left-invariant vector fields  
of Hodge type (1,0) on the corresponding nilpotent Lie group with 
respect to some  complex structure $I\in S^2$. A hypercomplex nilmanifold 
is called {\bf abelian} when this Lie algebra is abelian. We prove that 
all complex subvarieties of $(M,L)$ for generic $L\in S^2$ 
on a hypercomplex abelian nilmanifold are also hypercomplex nilmanifolds.
\end{minipage}
}

\tableofcontents


\section{Introduction}


\subsection{Complex geometry of non-algebraic manifolds}

Algebraic varieties are typically given by 
a system of polynomial equations. Complex non-algebraic
varieties are described in different and widely varying 
fashions, much less uniform.

One of the main sources of non-algebraic complex manifolds
comes from the theory of homogeneous or locally
homogeneous spaces. Many compact Lie groups, as well as
their homogeneous spaces, admit left-invariant complex 
structures (\cite{_Samelson_}, \cite{_Tits:homoge_},
\cite{_Joyce_}). When the group $G$ is not compact,
one can consider a cocompact lattice $\Gamma$ and take
(say) a left quotient $G/\Gamma$ equipped with a
left-invariant complex structure. Finally, one could
take a quotient of a domain $\Omega\subset \C^n$
by a holomorphic action of a discrete group.
This way one can obtain some of the non-K\"ahler 
complex surfaces (Hopf, Kodaira, Inoue), as well
as the multi-dimensional generalization of
the Inoue surfaces, called {\bf Oeljeklaus-Toma manifolds}
(\cite{_Oeljeklaus_Toma_}).

It is interesting that in most of these examples,
the complex geometry is more restrictive than
in the algebraic case. The Inoue surface
and some special classes of 
Oeljeklaus-Toma manifolds have no 
complex subvarieties at all
(\cite{_OV:Oeljeklaus-Toma_}).
More general Oeljeklaus-Toma manifolds
can admit subvarieties, but 
the subvarieties of smallest possible dimension
are also obtained as a quotient of a 
domain $\Omega\subset \C^n$ by an 
affine holomorphic action (\cite[Theorem 1.2]{_OVV:Oeljeklaus-Toma_flat_}). 
Moreover, the torsion-free flat connection,
intrinsic on the Oeljeklaus-Toma
manifolds, induces a torsion-free flat connection
on any of these subvarieties.

The complex structures on compact Lie groups and associated 
homogeneous spaces are also very interesting. By 
Borel-Remmert-Tits theorem, such manifolds
are always principal toric fibrations over
a rational projective manifold (\cite{_Tits:homoge_}, 
\cite{_Verbitsky:toric_}). It is possible to classify
all subvarieties of a class of principal toric fibrations,
called ``positive principal toric fibrations''.
In \cite{_Verbitsky:toric_} it was shown that
any irreducible subvariety of a positive principal toric fibration
is either a pullback from the base of the fibration, or
belongs to its fiber. 

An interesting example of a non-algebraic complex manifold
is a Bogomo\-lov--Guan manifold (\cite{_Guan_},
\cite{_Bogomolov_}). It is a non-K\"ahler holomorphic
symplectic manifold which is constructed from a Kodaira
surface in a similar way as generalized Kummer varieties
are constructed from complex $2$-dimensional tori (see
\cite{_Beauville_}, \cite{_Bogomolov_},
\cite{_KV:BG_}). A Bogomolov--Guan manifold is 
equipped with a natural Lagrangian
fibration over the projective space. A partial
classification of subvarieties of Bogomolov--Guan
manifold was obtained in
\cite{_Bogomolov_Kurnosov_Kuznetsova_Yasinsky_}, where the
algebraic dimension of Bogomolov--Guan manifolds was
computed as well.

It is natural to expect that these results can be
generalized. For example, consider a flat affine
complex manifold, that is, a manifold equipped 
with a flat torsion-free connection preserving
the complex structure. This class of manifolds
contains Hopf manifolds, Oeljeklaus-Toma manifolds
and compact tori. Subvarieties in a non-algebraic
torus and the Oeljeklaus-Toma manifolds are also
flat affine (\cite{_Ueno:albanese_}, 
\cite{_OV:Oeljeklaus-Toma_}), however, for
the Hopf manifold and the torus, such result
won't work, indeed, there are tori which are algebraic
and Hopf manifolds obtained as principal elliptic
bundles over a weighted projective space.

However, in both of these cases, it is enough
to deform the complex structure to a general point
to obtain a classification of subvarieties. As far as we know,
the following questions remains open. 

\hfill

\question
Let $M$ be a compact flat affine complex manifold.
Is it true that there exists a finite covering $\tilde M$ of $M$ 
and a flat affine complex manifold $M_1$ in the same
deformation class as $\tilde M$, such that all  complex
subvarieties of $M_1$ are also flat affine manifolds?

\hfill

A similar question can be stated in context of 
complex structures on locally homogeneous manifolds, defined in 
Subsection \ref{_complex_nilma_Subsection_}.

\hfill

\question\label{_subva_loh_hom_Question_}
Let $M= G/\Gamma$ be a compact locally homogeneous complex manifold.
Is it true that there exists a locally homogeneous complex manifold $M_1$ in the same
deformation class as $\tilde M$, such that all complex
subvarieties of $M_1$ are locally homogeneous 
submanifolds (\ref{_locally_homogeneous_submanifold_def_})?

\hfill

One of the easiest way to provide deformations of complex
structures is to use the quaternionic rotation. This is one of the
motivation of the present work, where we successfully 
solve \ref{_subva_loh_hom_Question_} 
in the presence of an abelian hypercomplex structure.

\subsection{Trianalytic subvarieties in hyperk\"ahler manifolds}

A great source of non-algebraic K\"ahler
manifolds is provided by the hyperk\"ahler geometry.
Consider a compact, K\"ahler, holomorphically
symplectic manifold $M$, such as a K3 surface.
Calabi-Yau theorem provides $M$ with a {\bf hyperk\"ahler
structure} (\ref{hyperkahler_def}), that is, a  triple 
$I,J,K$ of complex structure operators satisfying quaternionic
relations, and a Riemannian metric which is K\"ahler
with respect to $I,J,K$. One of the main tools of hyperk\"ahler
geometry is called ``the hyperk\"ahler rotation'':
any quaternion $L\in {\Bbb H}, L^2=-1$,
induces a complex structure on $M$. 
This gives the whole 2-dimensional sphere 
of complex strictures, called {\bf the twistor sphere},
or {\bf the twistor deformation} of a hyperk\"ahler
manifold.

\hfill

Starting from an algebraic manifold and applying the
rotation, we usually obtain non-algebraic ones.
Recall that the {\bf algebraic dimension} (Subsection 
\ref{algebraic_dimension_section}) $a(M)$
of a complex manifold $M$ is the transcendence
degree of its field of global meromorphic functions.
When $M$ is algebraic, one has $a(M)=\dim M$.
Hyperk\"ahler rotation provides examples
of complex manifolds with algebraic dimension 0:

\hfill

\proposition(\cite[Cor. 5.12]{_Fujiki_})
\label{algebraic_dimension_of_hyperkahler_prop}
Let $X$ be a compact hyperk\"ahler manifold. Then for all
but a countable number of complex structures $L\in\H$ the
algebraic dimension of the complex manifold $X_L$ is zero.

\hfill

This result is a part of a more general phenomenon,
which can be used to classify complex subvarieties
in general (non-algebraic) deformations of hyperk\"ahler
manifolds. For the further use, we formulate the following definition 
in a more general situation of hypercomplex geometry
(Subsection \ref{_hc_intro_Subsection_}).

\hfill

\definition
\label{trianalytic_def}
A subvariety $Z$ of a hypercomplex manifold $X$ is called
{\bf trianalytic} if it is complex analytic with respect
to every complex structure $L\in\H$.

\hfill

Trianalytic subvarieties enjoy many good properties:
they are completely geodesic 
(\cite{Verbitsky_deformations_of_trianalytic}), and can have only very
mild singularities, which are resolved by a normalization
(\cite{_Verbitsky:hypercomple_}). 

\hfill

\proposition
(\cite{Verbitsky_subvarieties_in_non_compact}, see also \cite{Verbitsky_subvarieties_in_compact})
\label{submanifolds_in_hyperkahler_prop}
Let $X$ be a hyperk\"ahler manifold (not necessarily
compact). Then for all but a countable number of complex
structures $L\in\H$, all compact complex subvarieties of
$X_L$ are trianalytic.

\hfill

This result was often applied in symplectic
geometry to infer that a general almost complex
structure on $X$ does not support complex curves
(\cite{_EV-balls_}, \cite{_Solomon_V:locality_}).

\hfill

\ref{algebraic_dimension_of_hyperkahler_prop} follows
directly from
\ref{submanifolds_in_hyperkahler_prop}. Indeed, let $X$ be
a compact hyperk\"ahler manifold. Let $L\in \H$ be a
complex structure such that $X_L$ does not contain
subvarieties which are not trianalytic. Then a fortiori
$X_L$ does not contain any divisors. Complex varieties
with no divisors are of algebraic dimension zero.

\subsection{Hypercomplex manifolds}
\label{_hc_intro_Subsection_}

Hyperk\"ahler manifolds are well understood, but there
are not many compact examples: so far, only two series and
two sporadic examples exist. It makes sense
to weaken the hyperk\"ahler condition and
consider {\bf hypercomplex manifolds}
(\ref{hypercomplex_nilmanifold_def}).
A hypercomplex structure is a triple of
complex structures satisfying quaternionic
indentities. Examples of hypercomplex
manifolds are abundant: there are homogeneous
examples (\cite{_Joyce_}), nilmanifolds and
solvmanifolds (\cite{_Barberis_Dotti_},
\cite{Dotti_Fino_hypercomplex}), fibered bundles
(\cite{_Verbitsky_HKT-exa_}), locally conformally hyperk\"ahler
examples (\cite{_Ornea:LCHK_}) and many others.

Hypercomplex manifolds were first introduced
by Ch. Boyer (\cite{_Boyer_}), who also gave the complete list
of compact hypercomplex manifold in real dimension 4:
 K3 surface, compact torus, and Hopf manifold.

\hfill

Any element $L\in {\Bbb H}$, $L^2=-1$ induces an
integrable complex structure on $M$, just like
in the hyperk\"ahler case. However,
\ref{algebraic_dimension_of_hyperkahler_prop} and
\ref{submanifolds_in_hyperkahler_prop} don't hold for
general hypercomplex manifolds. In fact, for
most examples listed above these results are false.
Indeed, by Borel-Remmert-Tits theorem, any 
complex manifold with trivial (or abelian) fundamental
group and admitting a transitive action of a 
compact Lie group is always a principal toric fibration over
a rational projective manifold (\cite{_Tits:homoge_}, 
\cite{_Verbitsky:toric_}). Therefore, any 
homogeneous hypercomplex manifold of this type which 
is not a torus always has positive algebraic
dimension. 

Another counterexample is given by the Hopf manifold
 $\H^n\setminus \{0\}/\Z$, where the $\Z$-action is given
by the multiplication by $\lambda\in\R^{>1}$. 
 For any complex structure $L\in\H$ the
manifold $X_L$ admits an elliptic fibration over $\C
P^{2n-1}$. Hence the algebraic dimension of $X_L$ is
always $2n-1$ and $X_L$ always contains an elliptic
curve. 

It is interesting that all counterexamples 
existing so far are hypercomplex manifolds with
non-trivial canonical bundle. In fact, the anticanonical
bundle for all these examples is globally generated
and semi-ample. 

We expect that a better control over subvarieties and
algebraic dimension of hypercomplex manifolds is possible 
when the canonical bundle is trivial.

Note that the canonical bundle of a hypercomplex manifold
can be interpreted in terms of its hypercomplex structure.
Recall that {\bf the Obata connection} of a hypercomplex
manifold $(M,I,J,K)$ is a torsion-free connection $\nabla$ 
such that $\nabla(I)=\nabla(J)=\nabla(K)=0$. Such a
connection exists and is unique on any hypercomplex
manifold (\cite{_Kaledin_},\cite{_Obata_}).
Clearly, the Obata connection acts
on the bundle $\Lambda^{2n,0}(M,I)$, which is 
equal to the canonical bundle of $(M,I)$ where $n=\dim_\H M$.
This action is compatible with the real structure
$\eta \arrow J(\bar \eta)$, hence its holonomy
is real. The corresponding character on the
holonomy group can be identified with the
determinant $\det:\; GL(n,\H) \arrow \R^{>0}$.
Whenever this character vanishes (or, equivalently,
when the holonomy of the Obata connection belongs to 
$\ker \det = SL(n,\H)$), the canonical class is also
trivial. 

\hfill

\conjecture\label{_SL_n_H_Conjecture_}
Let $(M,I,J,K)$ be a compact hypercomplex
manifold with the holonomy of Obata connection
in $SL(n,\H)$. Then for all complex structures 
$L\in {\Bbb H}$ outside of a countable subset,
the complex manifold $(M,L)$ has algebraic
dimension 0, and all its complex subvarieties are 
trianalytic.

\hfill

This suggestion is supported by a result of
\cite{Soldatenkov_Verbitsky_SL(n H)}, where 
trianalyticity was proven for codimension 1 and 2
complex subvarieties of a hypercomplex manifold
with holonomy in $SL(n,\H)$ in the presence
of an HKT-metric (\ref{HKT_def}). 

In \cite{Barberis_Dotti_Verbitsky}, it was shown
that any complex nilmanifold has trivial canonical bundle.
Then, it is natural to expect that 
\ref{_SL_n_H_Conjecture_} would hold for hypercomplex nilmanifolds.

This is the main purpose of the present paper.

\subsection{Hypercomplex nilmanifolds as iterated principal bundles}

Recall that {\bf a nilmanifold} (\ref{nilmanifold_def}) is a compact
manifold equipped with a transitive action of a connected
nilpotent Lie group. Equivalently (\cite{_Malcev_}), a nilmanifold
is a quotient $M=G/\Gamma$, where $G$ is a connected, 
simply connected nilpotent Lie group, and $\Gamma\subset G$
is a cocompact lattice (that is, a discrete subgroup
such that $G/\Gamma$ is compact). We shall usually consider
the left quotient; to emphasize this, we sometimes
write $M=\Gamma\backslash G$.

To give a left-invariant geometric structure (such as a complex
structure, or a hypercomplex structure) on a Lie group
is the same as to give it on its Lie algebra $\g =T_e G$
and extend to $G$ using the left action. In other words,
a left invariant almost complex structure on $G$
is the same as an operator $I:\; \g \arrow \g$
satisfying $I^2=-\Id$. Integrability of the
corresponding complex structure is  equivalent
to the relation 
\begin{equation}\label{_1,0_subalge_Equation_}
[\g^{1,0},\g^{1,0}]\subset \g^{1,0},
\end{equation}
where $\g^{1,0}=\{x\in \g_\R \otimes \C \ \ |\ \ I(x)=\1 x\}$.
Indeed, by Newlander-Nirenberg theorem, an almost complex
structure on a manifold is integrable whenever 
$[T^{1,0}M, T^{1,0}M]\subset T^{1,0}M$. However,
it suffices to check that $T^{1,0}M$ is generated
by vector fields $\xi_1, ..., \xi_n$ such that
$[\xi_i, \xi_j]\in T^{1,0}M$. This is obviously
implied by \eqref{_1,0_subalge_Equation_}.

\hfill

\example\label{_Kodaira_Example_}
As an example of a complex nilmanifold, we 
define the Kodaira surface, following
\cite{_Hasegawa_}. Consider the Lie algebra
spanned by vectors $x, y, z, t$, with
$[x,y]=z$ and all other commutators vanishing.
Let the complex structure act as
$I(x)=y$, $I(y)=-x$, $I(z)=t$ and $I(t)=-z$.
To show that the corresponding almost
complex structure is integrable, we
notice that the only possibly
non-trivial commutator between the
$(1,0)$-vectors, $[x+\1y, z+\1 t]$,
vanishes. From the exact sequence of Lie algebras
\[
0\arrow \langle t,z\rangle \arrow 
\langle x,y, t,z\rangle\arrow \langle x,y \rangle\arrow 0,
\]
it is possible to see that the Kodaira
surface admits a principal fibration over an elliptic
curve with elliptic fibers. The type of elliptic curves is determined
by the choice of the discrete lattice in the
corresponding Lie group, with the lattice different for
different Kodaira surfaces.

\hfill

By ``hypercomplex nilmanifold'' we always mean
a quotient  $M=\Gamma\backslash G$, where
$G$ is equipped with a left-invariant hypercomplex
structure. 

Let $Z\subset G$ be the center of a nilpotent Lie group
$G$, and $\Gamma\subset G$ a cocompact lattice. Denote 
the nilmanifold $\Gamma\backslash G$ by $M$. 
It is not very hard to see that $Z \cap \Gamma$ is
cocompact in $Z$, and $M_1:=\frac{G/Z}{\Gamma/(Z \cap \Gamma)}$
is a nilmanifold. Moreover, the natural projection
$M \arrow M_1$ is a fibration with the fiber
$\frac{Z}{Z \cap \Gamma}$, hence it is a principal
toric fibration. 

Repeating this construction, we prove that any
smooth nilmanifold can be obtained as an iterated
principal toric bundle, that is, a total space
of a principal toric fibration over 
a principal toric fibration over... ...over a point.

In complex category, this result is generally false.
A survey paper \cite{_Rollenske:iterated_toric_} 
by S. Rollenske  considers this question in some detail,
as well as his Ph. D. thesis \cite{_thesis:Rollenske_}.
Recall that the classification of 6-dimensional complex
nilmanifolds was obtained in \cite{_Abbena_Garbiero_Salamon_}.
In \cite[Example 1.14]{_thesis:Rollenske_}, Rollenske
produces an example of a 6-dimensional complex nilmanifold
which does not admit a structure of an iterated 
holomorphic principal bundle. In 
\cite[Corollary 3.10]{_Rollenske:iterated_toric_}
he proves that all other 6-dimensional complex 
nilmanifolds admit a structure of an iterated 
holomorphic principal bundle.

One of the most important observations made by Rollenske
deals with the abelian complex structures. Abelian complex 
structures were introduced by M. L. Barberis in her
dissertation \cite{_Barberis:PhD_}, and explored at some length by
Barberis, Dotti and their collaborators 
(\cite{_Barberis_Dotti_}, \cite{_ABD:classification_}, 
\cite{_Barberis_Dotti:solvable_}).
To appreciate their definition, note that the integrability
of a left-invariant complex structure operator is
equivalent to \eqref{_1,0_subalge_Equation_}, that
is, an almost complex structure $I\in \End(\g)$ is integrable
if and only if $\g^{1,0}$ is a Lie subalgebra in $\g \otimes_\R \C$. 
The complex structure $I$ is called {\bf abelian}
when the algebra $\g^{1,0}$ is abelian.

\hfill

S. Rollenske
(\cite[Cor. 3.11]{_Rollenske:iterated_toric_}) has shown
that all abelian complex nilmanifolds admit a structure of an iterated 
holomorphic principal bundle. In fact, their upper central
series are complex-invariant subspaces in $\g$ (\ref{upper_central_series_L_invariant_lemma}). 

\hfill

For a hypercomplex nilmanifold $(M,I,J,K)$,
the complex structure $I$  is abelian if and only
if all the complex structures induced by quaternions
are also abelian (\cite[Lemma 3.1]{Dotti_Fino_hypercomplex}). 
We show that the upper central series of an abelian
hypercomplex nilmanifold $M$ are ${\Bbb H}$-invariant
(\ref{upper_central_series_L_invariant_lemma}). 
This implies that $M$ is an iterated holomorphic principal
bundle, and, moreover, the corresponding projections
are holomorphic with respect to $I, J, K$.

\subsection{Subvarieties in hypercomplex nilmanifolds:
  main results}

In this subsection we state the main results of the present paper
and give a few hints about their proof. It is independent
from the main body of this paper. For a definition of 
a hypercomplex nilmanifold and more details, see
Subsection \ref{_hc_Subsection_}.

\hfill

The first of two main results of this paper:

\hfill

\theorem\label{_hc_alge_dim_Intro_Theorem_}
Let $(M,I,J,K)$ be a hypercomplex nilmanifold,
and 
\[ S:= \{L\in {\Bbb H}\ \ |\ \ L^2=-1\}
\] the 2-dimensional sphere of all complex
structures on $M$ induced by the quaternion action.
Then for all $L\in S$ outside of a countable
subset $R\subset S$, the algebraic dimension
of $(M,L)$ vanishes.

\proof \ref{algebraic_dimension_of_hypercomplex_cor}. \endproof

\hfill

The second result is based on
\ref{upper_central_series_L_invariant_lemma},
which gives a structural result about nilmanifolds
with abelian hypercomplex structure.

\hfill

\theorem\label{_hc_subvar_Intro_Theorem_}
Let $(M,I,J,K)$ be a nilmanifold
equipped with an abelian hypercomplex structure,
and  \[ S:= \{ L\in {\Bbb H}\ \ |\ \ L^2=-1\}
\] the 2-dimensional sphere of all complex
structures on $M$ induced by the quaternion action.
Then for all $L\in S$ outside of a countable
subset $R\subset S$, all complex subvarieties
$X\subset (M,L)$ are trianalytic and locally homogeneous.

\proof \ref{my_precious}. \endproof

\hfill

The proof of \ref{_hc_alge_dim_Intro_Theorem_}
is based on a result of \cite{GFV_algebraic_dimension}.
In this paper the authors consider a foliation
$\Sigma$ on a complex nilmanifold $M$ given by intersection of kernels
of all closed holomorphic 1-forms: 
\[ \Sigma=
\bigcap_{\eta\in \Lambda^{1,0}(M), d\eta=0} \ker \eta.
\]
The same foliation can be obtained by left-translates
of the subalgebra $[\g,\g]+ I ([\g,\g])$, where
$M=\Gamma\backslash G$, and $\g= \Lie(G)$. By 
\cite[Theorem 1.1]{GFV_algebraic_dimension},
any meromorphic function on $M$  is constant on the leaves of $\Sigma$.
However, the leaves of $\Sigma$ are not necessarily closed;
all meromorphic functions are constant on the leaves of the
foliation obtained by taking the closures of the leaves of $\Sigma$.

Generally speaking, the Lie algebra $\g$ of a 
nilmanifold $\Gamma\backslash G$ is equipped with a rational structure,
with the rational lattice generated by logarithms of $\gamma\in \Gamma$.
It is not hard to see that a foliation $\Sigma_V$ given by translates of a 
subspace $V\subset \g$ has closed leaves if and only if $V$ is rational.
Let $\hat V_\Q$ be the smallest rational subspace of $\g$
containing $V$. From what we stated, it is clear that 
$\Sigma_{\hat V_\Q}$ is a foliation with closed leaves,
and for each leaf of $\Sigma_V$, its closure is a leaf of
$\Sigma_{\hat V_\Q}$.

Let now $(M,L)$ be a complex nilmanifold, and $\Sigma_{L,\Q}$
the closed foliation obtained by translates of the smallest
$L$-invariant rational space $[\g,\g]_{L,\Q}$ containing $[\g,\g]+ L ([\g,\g])$.
From \cite[Corollary 1.11]{GFV_algebraic_dimension},
it follows that all meromorphic functions on $(M,L)$
are constant on the leaves of $\Sigma_{L,\Q}$,
and, moreover, factorize through the leaf space,
which is a complex torus.

Consider now a continuous family of complex structures
$L$ on $\g$ and the corresponding family of subspaces
$[\g,\g]_{L,\Q}$. Let $V$ be a vector space and 
$H(V)$ the set of all subspaces contained
in $V$. We call a map $\phi$ from a topological space to the set
of subspaces of a vector space $W$ 
{\bf lower semicontinuous} if the preimage of
$H(V)$ is closed for any vector space $V\subset W$.

Since the closure of a limit is
contained in the limit of closures, the map
$L\arrow [\g,\g]_{L,\Q}$ is lower semicontinuous.
However, there are only countably many possible
rational lattices in $\g$, and uncountably many
points in any continuous family of complex structures.
If we apply that argument to the twistor deformation
of a complex structure on a hypercomplex nilmanifold, 
we obtain that the function $L\arrow [\g,\g]_{L,\Q}$
is constant on an open subset of the set $S^2 \subset {\Bbb H}$
of induced complx structure. A more careful analysis of 
the quaternionic action implies that the function $L\arrow [\g,\g]_{L,\Q}$ is 
constant outside of a countable set (\ref{complex_rational_implies_H_invariant_lem}). 
Indeed, for generic $L$, the space $[\g,\g]_{L,\Q}$ 
is identified with the minimal rational quaternionic-invariant
subspace of $\g$ containing $[\g,\g]$.

Denote the corresponding foliation on $M=\Gamma\backslash G$ 
by $\Sigma_{{\Bbb H},\Q}$. The above argument implies that
all meromorphic functions on $(M,L)$ for any complex structure
$L\in {\Bbb H}$ outside of a countable set are constant on the 
leaves and factorize through the leaf space of $\Sigma_{{\Bbb H},\Q}$.
However, this leaf space is actually a hyperk\"ahler torus,
which has no global meromorphic functions for general $L$ by 
\ref{algebraic_dimension_of_hyperkahler_prop}.
This was a scheme of a proof of \ref{algebraic_dimension_of_hypercomplex_cor};
detailed argument is given in Subsection \ref{_Albanese_Subsection_}.

\hfill

To prove \ref{_hc_subvar_Intro_Theorem_},
we use the iterated fibration constructed in 
\ref{upper_central_series_L_invariant_lemma}.
From this lemma it follows that any
nilmanifold with an abelian hypercomplex structure
is represented as a principal toric bundle $\pi:\; M \arrow M_1$ over
a nilmanifold $M_1$ with an abelian hypercomplex structure,
and this fibration is trianalytic, that is,
holomorphic with respect to all complex structures
 induced by quaternions. 

Consider a subvariety $X\subset M$.
Using induction, we can assume that for all 
induced complex structures $L\in S^2\subset {\Bbb H}$, except a countable set,
the image $\pi(X)$ is trianalytic and locally homogeneous.
Replacing $M_1$ by $\pi(X)$, we may actually
assume that the restriction $\pi\restrict X$
is surjective. Applying \ref{submanifolds_in_hyperkahler_prop},
and excising more complex structures from $S^2$, we may assume that 
the fibers $T$ of the principal toric bundle $\pi$ are 
generic enough, and all irreducible complex subvarieties of these
fibers are trianalytic subtori (\ref{subvarieties_in_tori_prop}).

Then for a general $x\in M_1$ 
there exists a subtorus $T'\subset T$ such that
$\pi^{-1}(x)$ is a trianalytic subtorus, or a collection
of subtori, if $\pi^{-1}(x)$ has several irreducible components.
In the latter case, $T'$ can be considered as a function 
of a smooth point $y\in \pi^{-1}(x)$.
Clearly, the correspondence $y\mapsto T'$ is continuous
on the set of smooth points of $X$. Since $X$ is irreducible,
its smooth locus is connected, hence the torus $T'$
is actually constant as a function of $y$. Replacing $M$ by the quotient $M/T'$,
we can assume that $X$ is a meromorphic multisection of 
$\pi:\; M \arrow M_1$. The branch locus and the exceptional
locus of $\pi:\; X \arrow M_1$ produce divisors in $M_1$,
which is impossible, because by induction hypotesis,
all subvarieties of $M_1$ are trianalytic. Therefore,
the map $\pi:\; X \arrow M_1$ is etale. Passing to
a finite covering, we may assume that it is a section
of the principal bundle $\pi:\; M \arrow M_1$.
However, any section of a principal bundle
trivializes it, giving $M= M_1\times T$,
and its summand $X=M_1$ is trianalytic
because this decomposition is compatible
with a hypercomplex structure.

This was a scheme of a proof of  \ref{_hc_subvar_Intro_Theorem_};
detailed argument is given in Subsection
\ref{_Proof_Subsection_}.

The hypercomplex rotation makes the argument
much more accessible and easy to state. However, a similar
result can be proven for any iterated fibration of tori,
provided that these tori have algebraic dimension zero
and are sufficiently generic. The ``sufficienly
generic'' bit becomes ugly if one is interested
in the maximal generality, and we won't state
it here.


\section{Preliminaries}

\subsection{Complex nilmanifolds}
\label{_complex_nilma_Subsection_}




\definition Let $\g$ be a Lie algebra. 
\begin{itemize}
    \item The {\bf lower central series} of $\g$ is the descending filtration on $\g$
\begin{equation}
    ...\subset \g_i \subset ...\subset \g_1 \subset \g_0 = \g
\end{equation}
where $\forall i > 0\colon \g_i := [\g,\g_{i-1}]$. In particular, $\g_1 = [\g,\g]$ is the commutator subalgebra of $\g$.

    \item The {\bf upper central series} of $\g$ is the ascending filtration on $\g$
    \begin{equation}
        0 = \z^0 \subset \z^1 \subset ... \z^i \subset...
    \end{equation}
    where $\forall i\colon\z^i := \{x\in\g\:|\:[x,\g]\subset \z^{i-1}\}$. In particular, $\z^1 = \z(\g)$ is the center of $\g$.

\end{itemize}

If $G$ is a group one defines the {\bf lower and upper central series of $G$} in a similar way.

\hfill

\definition A Lie algebra $\g$ is called {\bf nilpotent} if the following holds: $\exists k\colon \g_k = 0$. The minimal number $k$ with this property is called the {\bf number of steps} of a nilpotent Lie algebra $\g$. If the number of steps of $\g$ equals $k$, then one calls $\g$ a {\bf $k$-step nilpotent Lie algebra}. In a similar fashion one defines a {\bf nilpotent group}, the {\bf number of steps of a nilpotent group} and a {\bf $k$-step nilpotent group}.

\hfill

The following proposition is well known. For the convenience of the reader, we give its proof here.

\hfill

\proposition
\label{lower_is_contained_in_upper_prop}
Let $\g$ be a $k$-step nilpotent Lie algebra, $0 = \g_k \subset \g_{k-1} \subset ...\g_1\subset \g_0 = \g$ be its lower central series and $0 = \z^0 \subset \z^1 ...$ be its upper central series. Then one has $\z^k = \g$ and $k$ is the minimal number with this property. Moreover, $\forall i = 0,...,k$ one has
\begin{equation}
    \g_{k-i} \subset \z^i
\end{equation}

\proof We prove the assertion by induction on the number $k$ of steps of $\g$. For $k=1$ (the case of abelian Lie algebras) the assertion is clear. Let $\g$ be a $k$-step nilpotent Lie algebra. Consider the quotient map $\pi\colon \g\to\g/\z^1 =: \tilde{\g}$. The Lie algebra $\tilde{\g}$ is $(k-1)$-step nilpotent. Define $\tilde{\z}^i$ and $\tilde{\g}_i$ to be respectively the upper and the lower central series of $\tilde{\g}$. One has
$$
\z^i = \pi^{-1}(\tilde{\z}^{i-1}), \quad \g_j\subset \pi^{-1}(\tilde{\g}_j)
$$
By induction hypothesis $\tilde{\z}^{k-1} = \tilde{\g}$ while $\tilde{\z}^{k-2}\varsubsetneqq\tilde{\g}$. Hence $\z^k = \pi^{-1}(\tilde{\z}^{k-1}) = \g$ while $\z^{k-1} = \pi^{-1}(\tilde{\z}^{k-2})\varsubsetneqq \g$, and the first assertion follows. By induction hypothesis one has $\tilde{\g}_{k-1-i}\subset \tilde{\z}^i$. We obtain that $\g_{k-i}\subset \pi^{-1}(\tilde{\g}_{k-i}) \subset \pi^{-1}(\tilde{\z}^{i-1}) = \z^i$.
\endproof{}

\hfill

\definition
\label{nilmanifold_def}
Let $G$ be a nilpotent Lie group and $\Gamma\subset G$ be a cocompact lattice\footnote{A {\bf cocompact lattice} $\Gamma$ is a discrete subgroup of $G$ such that the quotient $\Gamma\backslash G$ is compact.} in $G$. Then the quotient $X:=\Gamma\backslash G$ is called a {\bf nilmanifold.}

\hfill

It was shown by Mal'cev (\cite{_Malcev_}) that the nilpotent Lie group $G$ from \ref{nilmanifold_def} is uniquely detirmined by the nilmanifold $X$.

\hfill

\proposition(\cite{_Malcev_})
\label{Malcev_completion_prop}
Let $\Gamma$ be a nilpotent finitely generated group. Then there exists a simply connected nilpotent Lie group $G$ such that $\Gamma$ embeds into $G$ as a cocompact lattice. The group $G$ is unique up to isomorphism and depends functorially on $\Gamma$.

\hfill

\definition
\label{Malcev_completion_def} The group $G$ from \ref{Malcev_prop} associated to a nilpotent group $\Gamma$ is called the {\em Mal'cev completion of $\Gamma$}.

\hfill

\proposition(\cite{_Malcev_})
\label{Malcev_prop}
Let $X = \Gamma\backslash G$ be a nilmanifold, where $G$ is a simply connected nilpotent Lie group and $\Gamma\subset G$ a cocompact lattice. Then $\Gamma$ is naturally isomorphic to $\pi_1(X)$, and $G$ to the Mal'cev completion of $\pi_1(X)$.

\hfill

Let $\Gamma\subset G$ be a cocompact lattice. Consider the set $\Lambda := \{x\in\g \:|\: \exp(x)\in\Gamma\}$. One can show that $\Lambda$ is a maximal rank lattice in $\g$ (\cite{_Malcev_}). To be short, we will say that $\Lambda = \log(\Gamma)$ in the sequel.

\hfill

Let $L$ be an almost complex structure\footnote{An {\bf almost complex structure} $L$ on a vector space is a linear operator such that $L^2 = -1$.} on the Lie algebra $\g$ of $G$. Consider the decomposition $\g\otimes\C = \g^{1,0}\oplus\g^{0,1}$ into the direct sum of eigenspaces of $L$. The operator $L$ multiplies $\g^{1,0}$ by $\sqrt{-1}$ and multiplies $\g^{0,1}$ by $-\sqrt{-1}$. 

\hfill

\definition
An almost complex structure $L$ on $\g$ is called an {\bf
  integrable complex structure} or simply a {\bf complex
  structure} if $\g^{1,0}$ is a subalgebra of
$\g\otimes\C$.

\hfill

Be aware that the notion of a Lie algebra with an integrable complex structure $L$ is much broader than the notion of a complex Lie algebra. Indeed, a complex structure $L$ defines a structure of a complex Lie algebra on $\g$ if and only if
\begin{equation}
    \forall x,y\in\g\colon L[x,y] = [Lx,y] = [x,Ly],
\end{equation}
or equivalently, $[\g^{1,0},\g^{0,1}] = 0$ (then $\g^{1,0}$ is an ideal). If $L$ is an integrable complex structure on $\g$ then we denote the corresponding left-invariant complex structure on $G$ again by $L$. The group $G$ with the complex structure $L$ is, generally speaking, not a complex Lie group, because $L$ is not in general right-invariant. 

\hfill

\definition
\label{complex_nilmanifold_def}
Let $G$ be a nilpotent group with a left-invariant complex structure. Then the left quotient $X$ of $G$ by a cocompact lattice $\Gamma\subset X$ is called a {\bf complex nilmanifold}.

\hfill

\definition (\cite{_Barberis:PhD_}, see also \cite{_Barberis_Dotti_})
\label{abelian_complex_structure_def}
Let $\g$ be a Lie algebra with a complex structure $L$ and let $\g\otimes\C = \g^{1,0}\oplus\g^{0,1}$ be its eigenspace decomposition. If the subspace $\g^{1,0}$ is an abelian subalgebra of $\g\otimes\C$, then the complex structure $L$ is called an {\bf abelian complex structure} on $\g$.

\hfill

The subspace $\g^{1,0}$ consists of elements of the form $x-\sqrt{-1}Lx$, $x\in\g$. The abelianness of $L$ is equivalent to the equality
\begin{gather*}
    0 = [x-\sqrt{-1}Lx,y-\sqrt{-1}Ly] = \\
    =([x,y]-[Lx,Ly]) -\sqrt{-1}([Lx,y] + [x,Ly])
\end{gather*}
for every $x,y\in\g$. The real and the imaginary part of the expression above must both vanish. Therefore, a complex structure $L$ is abelian if and only if 

\begin{equation}
\label{_abelian_complex_structure_eq_}
    \forall x,y\in\g\colon [x,y] = [Lx,Ly]
\end{equation}
or equivalently,
\begin{equation}
\label{_abelian_complex_structure_eq_1_}
    \forall x,y\in\g\colon  [Lx,y] = -[x,Ly]
\end{equation}
\hfill

\definition
\label{_abelian_complex_nilmanifold_def}
Assume that the Lie algebra $\g$ is nilpotent. Equip the corresponding Lie group $G$ with the left-invariant complex structure induced by $L$. If $\Gamma\subset G$ is a cocompact lattice, then the nilmanifold $\Gamma\backslash G$ is called an {\bf abelian complex nilmanifold}.

\hfill

\definition
\label{hypercomplex_nilmanifold_def}
A {\bf hypercomplex structure} on a Lie algebra $\g$ is a triple of integrable complex structures $(I,J,K)$ on $\g$ satisfying quaternionic relations (see equation \eqref{quaternionic_relations_eq} below). We call a hypercomplex structure on a Lie algebra $\g$ {\bf abelian} if the complex structures $I,J,K$ are abelian. 

\hfill

By \cite[Lemma 3.1]{Dotti_Fino_hypercomplex} for any hypercomplex structure $(I,J,K)$ the complex structures $J$ and $K$ are abelian whenever $I$ is abelian. Therefore, in the definition above it is actually enough to require that only one of the complex structures $I,J,K$ is abelian.

The definition of {\bf hypercomplex nilmanifolds} and {\bf abelian hypercomplex nilmanifolds} is analogous to \ref{complex_nilmanifold_def}. The complex (resp. hypercomplex etc) structures on complex (resp. hypercomplex etc) nilmanifolds are examples of so called {\bf locally $G$-invariant structures}, as explained in the following definition.

\hfill

\definition Let $X = \Gamma\backslash G$ be a nilmanifold and $\g$ be the Lie algebra of $G$. Consider a tensor $\eta\in\g^{\otimes k}\otimes(\g^*)^{\otimes l}$. Let $\tilde{\eta}$ be the left-invariant tensor on $G$ induced by $\eta$. The tensor $\tilde{\eta}$ clearly descends to a tensor on $X$. In such a situation the induced tensor on a nilmanifold $X$ is called a {\bf locally $G$-invariant tensor} and the corresponding geometric structure is called a {\bf locally $G$-invariant geometric structure}.

\hfill

Be aware that a locally $G$-invariant tensor on a
nilmanifold $X = \Gamma\backslash G$ is generally not
invariant with respect to the right action of $G$ on $X$.
In other words, a complex nilmanifold is generally
speaking not a homogeneous complex manifold, though its
universal cover is $G$-homogeneous. We will call
such structures on $\Gamma\backslash G$ 
{\bf locally homogeneous complex structures.}


\subsection{Algebraic dimension}
\label{algebraic_dimension_section}

\definition
\label{algebraic_dimension_def}
Let $X$ be a compact complex manifold. Consider the field $K(X)$ of meromorphic functions on $X$. The transcendence degree of $K(X)$ over $\C$ is called the {\bf algebraic dimension} of $X$ and denoted by $\dim_{alg}(X)$

\hfill

One always has $\dim_{alg}(X)\le \dim(X)$. Recall that a complex manifold $X$ is called {\bf Moishezon} if $X$ is bimeromophic to a complex projective variety. The equality $\dim_{alg}(X) = \dim(X)$ holds if and only if $X$ is Moishezon (\cite{_Moishezon:1977_}).

Let $X$ be a compact complex manifold. By \cite{Campana_algebraic_reduction} there exists a meromorphic dominant map $r$ from $X$ to a projective variety $X^{alg}$ which induces an isomorphism $r^*\colon K(X^{alg})\to K(X)$ on the fields of meromorphic functions. The variety $X^{alg}$ can be constructed as follows. Let $(t_1,t_2,...t_N)$ be generators of $K(X)$ over $\C$ as a field. They define a meromorphic map $r\colon X\to \C P^N$ sending a point $x\in X$ outside of the poles of $t_i$'s to $[1:t_1(x):...:t_N(x)]\in\C P^N$. By \cite{Campana_algebraic_reduction} the closure $X^{alg}$ of the image of $r$ in $\C P^N$ is a projective variety such that the map $r^*\colon K(X^{alg})\to K(X)$ is an isomorphism.

\hfill

\definition 
\label{_algebraic_reduction_def_}
Let $X$ be a compact complex manifold. An {\bf algebraic reduction} of $X$ is a projective variety $X^{alg}$ together with a meromorphic dominant map $r\colon X\to X^{alg}$ such that the induced map $r^*\colon K(X^{alg})\to K(X)$ is an isomorphism.

\hfill

An algebraic reduction $X^{alg}$ of $X$ is defined uniquely up to a birational isomorphism. Moreover, every meromorphic map from $X$ to an algebraic variety is uniquely factorized through the meromorphic map $r\colon X\to X^{alg}$ (\cite{Campana_algebraic_reduction}).

\hfill

\lemma
\label{algebraic_dimension_finite_maps_lemma}
Let $f\colon X\to Y$ be a dominant map of complex manifolds. Then $\dim_{alg}(X) \ge \dim_{alg}(Y)$. If, moreover, the map $f$ is finite, then the equality holds.

\proof The map $f$ induces the embedding of the field $K(Y)$ into $K(X)$. Hence the transcedence degree of $K(Y)$ does not exceed that of $K(X)$. If $f$ is a finite map, then the extension $K(Y)\subset K(X)$ is finite. Hence the transcedence degrees of both fields coincide.
\endproof{}

\hfill

\subsection{Albanese variety}

Let $X$ be a compact complex manifold. Denote by $d\calO_X$ the sheaf of closed holomorphic differentials\footnote{A {\bf closed holomorphic differential} on a complex manifold $X$ is a closed $(1,0)$-form.} on $X$. Consider a short exact sequence of sheaves on $X$
\begin{equation}
    0\to \underline{\C}\to \calO_X \to d\calO_X\to 0
\end{equation}
It induces an embedding 
\begin{equation}
\label{_holomorphic_differentials_to_h1_eq_}
H^0(X,d\calO_X)\hookrightarrow H^1(X,\C)
\end{equation}
of cohomology groups. Let us dualize the map~\eqref{_holomorphic_differentials_to_h1_eq_} to get a surjection
\begin{equation}
\label{_h1_to_dual_holomorphic_differentials_eq_}
    H_1(X,\C)\twoheadrightarrow H^0(X,d\calO_X)^*
\end{equation}

\hfill

\definition(\cite{_Blanchard:albanese_},\cite{_Ueno:albanese_})
\label{_Albanese_def_}
Let $X$ be a compact complex manifold. Consider the composition of maps
\begin{equation}
\label{_albanese_eq_}
    H_1(X,\mathbb Z)\rightarrow H_1(X,\C) \twoheadrightarrow H^0(X,d\calO_X)^*
\end{equation}
The quotient $\Alb(X)$ of the space $H^0(X,d\calO_X)^*$ by
the minimal complex closed subgroup of $H^0(X,d\calO_X)^*$
containing the image of $H_1(X,\mathbb Z)$ is called the
{\bf Albanese variety} of $X$.

\hfill

\proposition
Let $X$ be a compact complex manifold. Then its Albanese
variety $\Alb(X)$ is a compact complex torus.

\hfill

\proof 
By definition, $\Alb(X)$ is a quotient of
$H^0(X,d\calO_X)^*$ by a closed subgroup $R$ which contains the
closure of the image of $H_1(X,\mathbb Z)$.  To prove that
$\Alb(X)$ is compact, it is enough to show that the
$\R$-vector space generated by $R$
in $H^0(X,d\calO_X)^*$ coincides with the whole vector
space $H^0(X,d\calO_X)^*$. Hence it suffices to
check that the
map~\eqref{_h1_to_dual_holomorphic_differentials_eq_}
sends $H_1(X,\R)$ surjectively onto $H^0(X,d\calO_X)^*$,
or equivalently, that the map
\begin{equation}
\label{_aa_}
    H^0(X,d\calO_X)\xrightarrow{\Re} H^1(X,\R)
\end{equation}
is injective. Consider a closed holomorphic differential
$\eta\in H^0(X,d\calO_X)$ such that $\Re \eta = df$ for
some smooth function $f$. Then $\eta = 2\partial f$ and
$\partial\bar{\partial}f = -\frac{1}{2}\bar{\partial}\eta
= 0$. We obtain that $f$ is a pluri-harmonic function on a
compact complex manifold, hence constant by maximum
principle. Hence $\eta = 2\partial f = 0$. We have proved
that the map~\eqref{_aa_} is injective. \endproof{}

\hfill

Fix a point $x_0\in X$. For any point $x\in X$ choose a path $\gamma$ connecting the points $x_0$ and $x$. The path $\gamma$ defines a functional $a_\gamma$ on the space $H^0(X,d\calO_X)$ of closed holomorphic differentials on $X$ in the following way:
\begin{equation}
    a_\gamma (\eta) := \int\limits_{\gamma} \eta
\end{equation}

If $\gamma'$ is another path connecting the points $x_0$ and $x$ then the functional $a_{\gamma'}-a_\gamma\in H^0(X,d\calO_X)^*$ lies in the image of $H_1(X,\mathbb Z)$ in $H^0(X,d\calO_X)^*$ by the map~\eqref{_albanese_eq_}. We obtain a correctly defined map
\begin{equation}
\label{_albanese_map_eq_}
    A\colon X\to \Alb(X)\quad x\mapsto a_\gamma
\end{equation}
where $\gamma$ is any path connecting $x_0$ and $x$.

\hfill

\definition
\label{_albanese_map_def_}
Let $X$ be a compact complex manifold. The map $A\colon X\to \Alb(X)$ defined in \eqref{_albanese_map_eq_} is called the {\bf Albanese map} of $X$.

\hfill

One can show that the Albanese map $A\colon X\to \Alb(X)$ is holomorphic (\cite{_Blanchard:albanese_}).

Let $X = \Gamma\backslash G$ be a nilmanifold
(\ref{complex_nilmanifold_def}).
Consider a left-invariant vector field $\tilde{\xi}$  on $G$
corresponding to a vector $\xi\in\g$.
Then $\tilde{\xi}$ descends to a locally $G$-invariant vector field on $X = \Gamma\backslash G$ denoted also by $\tilde{\xi}$. Let $\Vol_X$ be the volume form on $X$ induced by a left-invariant volume form on $G$ such that $\int_X\Vol_X=1$. Given a $1$-form $\eta$ on $X$ we define a covector $\eta_{\inv}\in\g^*$ as follows
\begin{equation}
\label{_invariant_form_eq_}
    \eta_{\inv}(\xi) := \int\limits_X \eta(\tilde{\xi})\Vol_X
\end{equation}

\hfill

One can show that if $\eta$ is a closed 
$1$-form on $X$ then $\eta_{\inv}$ is a closed $1$-form on
$\g$. Moreover, if $\eta$ is exact, the form 
$\eta_{\inv}$ vanishes  (\cite{_Nomizu_}).
Therefore, the map sending a form $\eta\in\Gamma(\Omega^1X)$ to $\eta_{\inv}\in\g^*$ descends to a map $\operatorname{Av}\colon H^1(X,\R)\to H^1(\g,\R)$ of cohomology groups. 

\hfill

\definition (\cite{GFV_algebraic_dimension})
Let $X = \Gamma\backslash G$ be a nilmanifold. The map 

\begin{equation}
    \operatorname{Av}\colon H^1(X,\R)\to H^1(\g,\R)= \left(\frac{\g}{[\g,\g]}\right)^*
\end{equation}
sending the cohomology class of a $1$-form $\eta$ on $X$ to the cohomology class of the form $\eta_{\inv}\in\g^*$ is called the {\bf averaging map}.

\hfill

\proposition (\cite{_Nomizu_}) Let $X=\Gamma\backslash G$ be a nilmanifold. Every closed $1$-form $\eta$ on $X$ is cohomologous to the locally $G$-invariant $1$-form on $X$ associated to $\eta_{\inv}\in\g^*$. Moreover, the averaging map $\operatorname{Av}\colon H^1(X,\R)\to H^1(\g,\R)$ is an isomorphism. \endproof{}

\hfill

Let $X = \Gamma\backslash G$ be a complex nilmanifold. Denote by $\Lambda =\log(\Gamma)\subset \g$ a maximal rank lattice in $\g$ and by $L$ the complex structure on $\g$. For any $\R$-subspace $W\subset \g$ let us denote by $W_{\Q,L}$ the minimal rational $L$-invariant subspace of $\g$ containing $W$. Here we call a subspace $W\subset V$ {\bf rational} if $\Lambda\cap W$ is a maximal rank lattice in $W$. The subspace $[\g,\g]_{\Q,L}$ is a normal subalgebra of $\g$, indeed, one has the following inclusions:
$$
[\g,[\g,\g]_{\Q,L}]\subset [\g,\g]\subset [\g,\g]_{\Q,L}
$$
The quotient Lie algebra $\mathfrak t:=\g/[\g,\g]_{\Q,L}$ is an abelian Lie algebra. Let $\Lambda_T$ denote the image of $\Lambda$ in $\mathfrak t$. The subgroup $\Lambda_T\subset \t$ is a lattice of maximal rank in $\t$. The quotient map $p\colon \g\to \mathfrak t$ induces a map on nilmanifolds:
\begin{equation}
\label{Albanese_map_eq}
    P\colon X\to T:=\mathfrak t/\Lambda_T
\end{equation}

\hfill

\proposition (see also \cite{Rollenske_albanese})
\label{is_albanese_prop}
Let $X = \Gamma\backslash G$ be a complex nilmanifold. The complex torus $T$ defined in equation (\ref{Albanese_map_eq}) is the Albanese variety of $X$ (\ref{_Albanese_def_}).

\hfill

\proof The averaging map $\operatorname{Av}\colon H^1(X,\C)\to H^1(\g,\C)$ sends the subspace $H^0(X, d\calO_X)\subset H^1(X,\C)$ isomorphically onto the subspace
$$
\left\{\alpha\in (\g^{1,0})^*\mid d\alpha = 0\right\} = \left(\frac{\g}{[\g,\g]+L[\g,\g]}\right)^* \subset H^1(X,\C)
$$

The image of $H_1(X,\mathbb Z)$ in the space
$H^0(X,d\calO_X)^* = \g/([\g,\g]+L[\g,\g])$ by the
map~\eqref{_albanese_eq_} is precisely the image of the
lattice $\Lambda\subset \g$ in this space
(\cite{_Nomizu_}). Consider a complex subgroup
$[\g,\g]_{\Q,L} + \Lambda$ of $\g$. This subgroup is
mapped onto a closed complex subgroup of
${\g/([\g,\g]+L[\g,\g])}$ under the quotient map ${\g \to
  \g/([\g,\g] + L[\g,\g])}$. Moreover, the image of
$[\g,\g]_{\Q,L} + \Lambda$ in $H^0(X,d\calO_X)^*$ is the
minimal closed complex subgroup of $H^0(X,d\calO_X)^*$
containing the image of $H_1(X,\mathbb Z)$. Hence,
\begin{equation}
    \Alb(X) = \frac{\g}{[\g,\g]_{\Q,L}+\Lambda}
\end{equation}
as desired. \endproof{}

\hfill



\proposition (\cite{GFV_algebraic_dimension})
\label{maps_to_kahler_prop}
Every meromorphic map
from a complex nilmanifold $X$ to a K\"ahler manifold is uniquely factorized
through the Albanese map (\ref{Albanese_map_eq}).

\hfill

As an immediate corollary of \ref{maps_to_kahler_prop} one has:

\hfill

\corollary (\cite{GFV_algebraic_dimension})
\label{algebraic_dimension_of_complex_nilmanifolds_cor}
Let $X$ be a complex nilmanifold and $T$ the Albanese variety of $X$. Then the following holds.
\begin{equation}
    \dim_{alg}(X) = \dim_{alg}(T)
\end{equation}

\proof \ref{algebraic_dimension_finite_maps_lemma} implies that $\dim_{alg}(X) \ge \dim_{alg}(T)$. Let $r\colon X \to X^{alg}$ be an algebraic reduction of $X$. The map $r$ factorizes uniquely through the Albanese map $P\colon X\to T$ as follows from \ref{maps_to_kahler_prop}. We use \ref{algebraic_dimension_finite_maps_lemma} again to obtain that $\dim_{alg}(T)\ge \dim_{alg}(X^{alg}) = \dim_{alg}(X)$.
\endproof{}

\hfill


\subsection{Hypercomplex manifolds}
\label{_hc_Subsection_}


Recall that the algebra $\H$ of quaternions is a $4$-dimensional algebra generated over $\R$ by $I,J,K$ subject to the relations
\begin{equation}
\label{quaternionic_relations_eq}
I^2 = J^2 = K^2 = -1;\:\:\: IJ = -JI = K
\end{equation}
An element $L = t + aI + bJ + cK \in H$ satisfies $L^2 = -1$ if and only if $t=0$ and $a^2 + b^2 + c^2 = 1$. Therefore such elements form a $2$-dimensional sphere inside of $\H$.

\hfill

\definition
Let $X$ be a smooth manifold equipped with a linear action of the algebra $\H$ on $TX$ i.e. an $\R$-linear map $\H\to \End(TX)$. Such a manifold is called {\bf almost hypercomplex}. Every element $L\in\H$ such that $L^2=-1$ induces an almost complex structure on $X$ which we denote by the same symbol $L$. If every complex structure $L$ on $X$ induced from quaternions is integrable, $X$ is called a {\bf hypercomplex manifold}. 

\hfill

In the definition of a hypercomplex manifold it is
actually enough to require the integrability of only two
linearly independent complex structures
(\cite{_Kaledin_},\cite{_Obata_}). For any complex
structure $L\in\H$ we denote by $X_L$ the corresponding
complex manifold.

Every (almost) hypercomplex manifold $X$ admits a {\bf hyper-Hermitian metric} $g$ i.e. a Riemannian metric which is Hermitian with respect to $I,J,K$. One can obtain such a metric by averaging an arbitrary metric on $X$ by the group $SU(2)$ of quaternions of unit norm. 

\hfill

\definition
A hypercomplex manifold equipped with a fixed hyper-Hermitian metric is called a {\bf hyper-Hermitian manifold}. 

\hfill 

For every complex structure $L\in \H$ on a hyper-Hermitian manifold $X$ one defines a $2$-form $\omega_L$ on $X$ as $\omega_L(x,y) := g(Lx,y)$. 

\hfill

\definition
\label{hyperkahler_def}
A hyper-Hermitian manifold $X$ is called {\bf hyperk\"ahler} if the $2$-forms $\omega_L$ are closed for every complex structure $L\in\H$.

\hfill 

Define $\Omega_I:=\omega_J+\sqrt{-1}\omega_K$. One can easily check that $\Omega_I$ is a $(2,0)$-form. 

\hfill 

\definition 
\label{HKT_def}
A hyper-Hermitian manifold $X$ is called an {\bf HKT-manifold} ("hyperk\"ahler with torsion") if $\partial \Omega_I = 0$. Here $\partial\colon \Lambda^{p,q}_I X \to \Lambda^{p+1,q}_I X$ is a $(1,0)$-part of the de Rham differential.

\hfill

Be aware that the HKT-condition (\ref{HKT_def}) is weaker than the hyperk\"ahler condition (\ref{hyperkahler_def}). 

\hfill

A complex nilmanifold $X = \Gamma\backslash G$ admits a K\"ahler structure if and only if $X$ is a complex torus (\cite{BG_kahler_nilpotent}). Hence the only example of a nilmanifold with a hyperk\"ahler structure is a hypercomplex torus, i.e. a quotient of a quaternionic vector space by a maximal rank lattice. It was shown in \cite{Fino_Grantcharov} that a hypercomplex nilmanifold $X = \Gamma\backslash G$ admits an HKT-metric if and only if it admits a $G$-invariant HKT-metric. 

\hfill

\proposition(\cite[Thm. 4.6]{Barberis_Dotti_Verbitsky}, \cite[Prop. 2.1]{Dotti_Fino_abelian_hypercomplex_are_hkt})
\label{hypecomplex_abelian_equals_HKT}
Let $X = \Gamma\backslash G$ be a hypercomplex nilmanifold with an invariant hyper-Hermitian metric $g$. The metric $g$ is HKT if and only if the hypercomplex structure on $X$ is abelian.

\hfill


\section{$\H$-Albanese variety}
\subsection{Definition of an $\H$-Albanese variety}

For a compact hypercomplex manifold $X$ consider the vector space $(\Omega^1X)_{par}$ of real $1$-forms ${\eta\in\Omega^1X}$ such that
\begin{equation}
\label{_closed_eq_}
    d\eta = dI\eta = dJ\eta = dK\eta = 0
\end{equation}

\hfill

\proposition
\label{_obata_par_prop_}
Let $X$ be a hypercomplex manifold. A $1$-form $\eta$
satisfies the condition~\eqref{_closed_eq_} if and only if
$\eta$ is parallel with respect to the Obata connection
$\nabla$ on $X$. \footnote{For the definition of Obata
  connection, see Subsection \ref{_hc_intro_Subsection_}.}

\hfill

\proof For any tensor $\alpha\in\Omega^1\otimes\Omega^1 =
S^2(\Omega^1X)\oplus \Omega^2X$ consider the decomposition
$\alpha = \alpha_{sym} + \alpha_{skew}$ where
$\alpha_{sym}\in S^2(\Omega^1X)$ and
$\alpha_{skew}\in\Omega^2X$. The Obata connection $\nabla$
is torsion-free, hence for any $1$-form $\eta$ one has
$(\nabla\eta)_{skew} = \frac{1}{2}d\eta$. 

Assume that $\eta$ is a $\nabla$-parallel $1$-form. Then
$d\eta = 2(\nabla\eta)_{skew} = 0$. For any complex
structure $L\in\H$ one has $\nabla L\eta = L\nabla\eta =
0$. Hence $dL\eta = 0$.

Assume that $\eta$ satisfies the
condition~\eqref{_closed_eq_}. Let $\eta^{1,0}$ be the
$(1,0)$-part of $\eta$ with respect to the complex
structure $I$. By \cite[Prop. 2.2]{_Soldatenkov:holonomy_} we have
\begin{equation}
    \nabla \eta^{1,0} = \bar{\partial}\eta^{1,0} - J\partial J\eta^{1,0} = 0
\end{equation}
Hence $\nabla \eta = 2\Re(\nabla\eta^{1,0}) = 0$. The
proposition is proved. \endproof{}

\hfill

\ref{_obata_par_prop_} justifies the notation $(\Omega^1X)_{par}$ introduced above.

The space $(\Omega^1X)_{par}$ is preserved by the linear
action of $\H$ on $1$-forms. Take a $1$-form $\eta\in
(\Omega^1X)_{par}$. For every complex structure $L\in\H$
the form $\eta$ is the real part of a closed
$L$-holomorphic $1$-form $\eta^{1,0}_L:=\frac{\eta +\1  L\eta}{2}$. 
Denote by $\calO^L_X$ the sheaf
of holomorphic functions on the complex manifold
$X_L$. For every complex structure $L\in\H$ we have the
following embeddings
\begin{equation}
    (\Omega^1X)_{par}\xlongrightarrow{\eta\ \mapsto \ \eta^{1,0}_L} H^0(X,d\calO^L_X) \stackrel{\Re}{\hookrightarrow} H^1(X,\R)
\end{equation}

\hfill

We introduce the following definition.

\hfill

\definition
\label{_h_albanese_def_}
Let $X$ be a compact hypercomplex manifold. Consider the composition of maps
\begin{equation}
    H_1(X,\mathbb Z)\arrow H_1(X,\mathbb R)\twoheadrightarrow ((\Omega^1X)_{par})^*
\end{equation}

The quotient $\Alb_\H(X)$ of the space
$((\Omega^1X)_{par})^*$ by the minimal trianalytic closed
subgroup of $((\Omega^1X)_{par})^*$ containing the image
of $H_1(X,\mathbb Z)$ is called the {\bf $\H$-Albanese
  variety} of $X$. 

\hfill


\hfill

For any complex structure $L\in\H$ the surjection $H^0(X,d\calO^L_X)^*\twoheadrightarrow ((\Omega^1X)_{par})^*$ induces a holomorphic map

\begin{equation}
    \phi_L\colon \Alb(X_L) \arrow (\Alb_\H(X))_L
\end{equation}

Fix a point $x_0\in X$. Similarly to the complex case we can construct a map

\begin{equation}
    A_\H\colon X\to \Alb_\H(X)
\end{equation}

The map $A_\H$ is defined as follows. For any $x\in X$ choose a path $\gamma$ connecting $x_0$ and $x$. The path $\gamma$ defines a functional $a_\gamma\in ((\Omega^1X)_{par})^*$ as
$$
a_\gamma(\eta) = \int\limits_{\gamma}\eta\quad \forall\eta\in(\Omega^1X)_{par}
$$
The map $A_\H$ sending $x\in X$ to $a_\gamma$ is a
well-defined map from $X$ to $\Alb_\H(X)$. This map is
called the {\bf $\H$-Albanese map} of $X$.

By construction, for every complex structure $L\in\H$ the $\H$-Albanese map $A_\H\colon X\to \Alb_\H(X)$ is equal to the composition of maps
\begin{equation}
    X \xrightarrow{A_L} \Alb(X_L) \xrightarrow{\phi_L} \Alb_\H(X)
\end{equation}
where $A_L$ is the Albanese map of $X_L$. Both maps $A_L\colon X\to \Alb(X_L)$ and $\phi_L\colon \Alb(X_L)\to (\Alb_\H(X))_L$ are holomorphic. It follows that the $\H$-Albanese map preserves the hypercomplex structure.


\subsection{Albanese varieties of a hypercomplex nilmanifold}
\label{_Albanese_Subsection_}

Let $V$ be a real vector space equipped with an action of $\H$ by linear endomorphisms and $\Lambda\subset V$ be a lattice of maximal rank. Recall that a real subspace $W\subset V$ is called {\bf rational} if $\Lambda\cap W$ is a maximal rank lattice in $W$. Let $L\in\H$ be a complex structure. For any real vector subspace $W\subset V$ let $W_{\Q,L}$ denote the minimal rational $L$-invariant subspace of $V$ containing $W$ and let $W_{\Q,\H}$ denote the minimal rational $\H$-invariant subspace of $V$ containing $W$ (see also Subsection~\ref{algebraic_dimension_section}).

\hfill

\lemma
\label{complex_rational_implies_H_invariant_lem}
Let $V$ be a real vector space with a linear action of $\H$, $\Lambda\subset V$ a maximal rank lattice. Then for all but at most countable number of complex structures $L\in\H$ the following property is satisfied: if $W\subset V$ is a rational $L$-invariant subspace of $V$ then $W$ is $\H$-invariant.

\hfill

\proof Let $W\subset V$ be a subspace of $V$. If $W$ is invariant with respect to two linearly independent complex structures $L,L'\subset V$ then $W$ is $\H$-invariant. Indeed, every two linearly independent complex structures generate $\H$ as an algebra. Therefore, every subspace $W\subset V$ is either not invariant with respect to all complex structures $L\in\H$, or invariant with respect to exactly two complex structures $L,-L\in\H$, or $\H$-invariant. The set of rational subspaces $W$ which are invariant with respect to exactly two complex structures $L_W,-L_W\in\H$ is countable. Let us exclude the complex structures which arise as such $L_W$. For every other complex structure $L\in\H$ each rational $L$-invariant subspace is $\H$-invariant.
\endproof{}

\hfill

\lemma
\label{lemma_on_rational_subspaces}
Let $W\subset V$ be a real subspace, $W_{\Q,L}$ and $W_{\Q,\H}$ be as above. Then for all but at most countable number of complex structures $L\in\H$ one has $W_{\Q,L} = W_{\Q,\H}$.

\hfill

\proof
By \ref{complex_rational_implies_H_invariant_lem} for all but at most countable number of complex structures $L\in\H$ the vector subspace $W_{\Q,L}$ is $\H$-invariant. If this is the case then one has $W_{\Q,L} = W_{\Q,\H}$.
\endproof{}

\hfill

Let $X = \Gamma\backslash G$ be a hypercomplex nilmanifold (\ref{hypercomplex_nilmanifold_def}), $\g$ the Lie algebra of $G$, and $\Lambda:=\log(\Gamma)$ the lattice in $\g$. We consider the abelian Lie algebra $\mathfrak t:= \g/[\g,\g]_{\Q,\H}$ with a hypercomplex structure (compare with Subsection~\ref{algebraic_dimension_section} where $\mathfrak t$ was defined as $\g/[\g,\g]_{\Q,L})$. Let $\Lambda_T$ be the image of the lattice $\Lambda$ in $\t$. The quotient map $p\colon \g\to\t$ induces a map 
\begin{equation}
\label{_a_}
P_\H\colon X\to T:=\t/\Lambda_T
\end{equation}

It follows from the construction that the map $P_\H$ is compatible with the hypercomplex structure.

\hfill

\proposition Let $X = \Gamma\backslash G$ be a hypercomplex nilmanifold. The hypercomplex torus $T$ defined in \eqref{_a_} is the $\H$-Albanese variety of $X$ (\ref{_h_albanese_def_}).

\proof The proof follows the same lines as the proof of \ref{is_albanese_prop}.
\endproof

\hfill

\theorem 
\label{albanese_of_hypercomplex_thm}
Let $X$ be a hypercomplex nilmanifold and $T$ the $\H$-Albanese variety of $X$. Then for all but at most countable number of complex structures $L\in\H$ the map $P_\H\colon X\to T$, considered as a morphism of complex manifolds $P_\H\colon X_L\to T_L$, is the Albanese map of $X_L$.

\hfill

\proof
\ref{lemma_on_rational_subspaces} applied to the subspace $[\g,\g]\subset \g$ implies that for all but at most countable number of complex structures $L\in\H$ one has $[\g,\g]_{\Q,L} = [\g,\g]_{\Q,\H}$. By the construction in Section \ref{algebraic_dimension_section} the Albanese variety of $X_L$ is the quotient of $\g/[\g,\g]_{\Q,L}$ by a lattice.
\endproof{}

\hfill

\conjecture \ref{albanese_of_hypercomplex_thm} is true for any hypercomplex manifold $X$.


\hfill



The following statement is an analogue of \ref{algebraic_dimension_of_hyperkahler_prop} for hypercomplex nilmanifolds.

\hfill

\corollary
\label{algebraic_dimension_of_hypercomplex_cor}
Let $X = \Gamma\backslash G$ be a hypercomplex nilmanifold. Then for all but at most countable complex structures $L\in\H$ one has $\dim_{alg}(X_L) = 0$.

\hfill

\proof By
\ref{algebraic_dimension_of_complex_nilmanifolds_cor} the
algebraic dimension of $X_L$ is equal to the algebraic
dimension of the Albanese variety of $X_L$. By
\ref{albanese_of_hypercomplex_thm} for all but at most
countable number of complex structures $L\in\H$ the
Albanese variety of $X_L$ is isomorphic to $T_L$ where $T$
is the  ${\Bbb H}$-Albanese variety of $X$. Any
translation-invariant hyper-Hermitian metric on a
hypercomplex torus is
hyperk\"ahler. \ref{algebraic_dimension_of_hyperkahler_prop}
implies $\dim_{alg}(T_L) = 0$ for all but at most
countable number of complex structures
$L\in\H$. \endproof{}


\section{Abelian (hyper-)complex varieties and their Albanese varieties}

The aim of this section is to deduce several useful observations about the geometry of abelian hyper(complex) varieties from a collection of linear algebra results. As before, we denote by $\g$ a nilpotent Lie algebra, by $\z^i$ the elements of its upper central series and by $\g_i$ the elements of its lower central series.

\hfill

\lemma
\label{upper_central_series_L_invariant_lemma}
Let $\g$ be a nilpotent Lie algebra with an abelian complex structure $L$ (\ref{abelian_complex_structure_def}).
Consider the upper central series $0=\z^0\subset\z^1\subset...\subset\z^k = \g$ of $\g$. Then for all $i=0,...k$, the subalgebras $\z^i$ are $L$-invariant. If $\g$ is equipped with an abelian hypercomplex structure $(I,J,K)$ then the subalgebras $\z^i$ are $\H$-invariant.

\proof
The proof goes by induction on $i$. The statement trivially holds for $i=0$. Assume we know the statement for the subalgebra $\z^i$. Take an element $x\in \z^{i+1}$. For any $y\in\g$ one has $[Lx,y] = -[x,Ly]$ as the complex structure $L$ is abelian (see equation~(\ref{_abelian_complex_structure_eq_1_})). As $x$ is contained in $\z^{i+1}$, one has $-[x,Ly]\in \z^i$. Therefore, $Lx$ also lies in $\z^{i+1}$. Let us assume that $\g$ is equipped with an abelian hypercomplex structure. For every $i$ the subalgebra $\z^i$ is $L$-invariant for each complex structure $L \in\{I,J,K\}$. Hence $\z^i$ is $\H$-invariant. \endproof{}

\hfill

Be aware that the statement of the lemma, generally
speaking, does not hold for lower central series. Indeed,
an element of the form $L[x,y]$ does not necessarily lie
in the commutator of $\g$, even for a Kodaira surface. 
Indeed, consider the Kodaira surface, determined
by the Lie algebra $\g \langle x, y, z, t\rangle$
as in \ref{_Kodaira_Example_}. Then 
$[\g, \g]=\langle z\rangle$, and $I(z)\notin [\g,\g]$.

\hfill

\lemma
\label{central_series_rational_lemma}
Let $\g$ be a nilpotent Lie algebra equipped with a maximal rank lattice $\Lambda\subset \g$ invariant with respect to the Lie bracket. Then the terms $\g_i$ of the lower central series of $\g$ as well as the terms $\z^i$ of the upper central series of $\g$ are rational subspaces.

\proof Clear, because $\g$ is a rational algebra. \endproof{}

\hfill

\lemma
\label{_principal_toric_fiber_bundle_lemma_}
Let $X = \Gamma\backslash G$ be a complex nilmanifold, $\g$ the Lie algebra of $G$ and $\z\subset \g$ the center of $\g$. Consider a rational $L$-invariant central subalgebra $\z'\subset \z$ and define $Z':=\exp(\z')$. Then the holomorphic map
\begin{equation}
    \pi\colon X = \Gamma\backslash G \arrow  Y:= \frac{G/Z'}{\Gamma/(\Gamma\cap Z')}
\end{equation}
is a holomorphic principal toric fiber bundle. Its fiber is $T':=Z'/(\Gamma\cap Z')$.

\hfill

\proof We claim that the right action of $Z'$ on $G$ is holomorphic i.e. the map
\begin{equation}
    R\colon G\times Z' \to G\quad (g,z)\mapsto gz
\end{equation}
is holomorphic. Indeed, $R$ is holomorphic in $z\in Z'$ because the left action of $G$ on itself preserves the complex structure. The map $L$ is holomorphic in $g\in G$ because the right action of $Z'$ on $G$ coincides with the left action. The holomorphic right action of $Z'$ on $G$ descends to a holomorphic action of $Z'/(\Gamma\cap Z')$ on ${Y = \frac{G/Z'}{\Gamma/(\Gamma\cap Z')}}$. The group $Z'/(\Gamma\cap Z')$ acts transitively on the fibers of $\pi\colon X\to Y$, hence $\pi$ is a principal toric fiber bundle. \endproof{}

\hfill

\proposition
\label{albanese_non_trivial_cor}
Let $X = \Gamma\backslash G$ be an abelian complex nilmanifold (\ref{abelian_complex_structure_def}). Then 
\begin{description}
    \item[(i)] The Albanese variety of $X$ is positive dimensional.
    \item[(ii)](\cite[Cor. 3.11]{_Rollenske:iterated_toric_}) 
There exists a finite sequence of holomorphic maps
    \begin{equation}
        X = X_0 \to X_1 \to...\to X_{N-1}\to X_N = \{pt\}
    \end{equation}
    such that for every $i$ the manifold $X_i$ is an abelian complex nilmanifold (\ref{_abelian_complex_nilmanifold_def}) and the map $X_i\to X_{i+1}$ is a principal toric fiber bundle.
\end{description}
{\bf Proof of (i):} By \ref{lower_is_contained_in_upper_prop} the commutator subalgebra $[\g,\g]$ is contained in $\z^{k-1}$ where $k$ is the number of steps of $\g$. \ref{upper_central_series_L_invariant_lemma} together with \ref{central_series_rational_lemma} imply that $\z^{k-1}$ is a proper $L$-invariant rational subspace of $\g$. Hence the minimal rational $L$-invariant subspace containing $[\g,\g]$, which we denote by $[\g,\g]_{\Q,L}$ as before, is contained in $\z^{k-1}$. In particular, $[\g,\g]_{\Q,L}$ is a proper subspace of $\g$. If $\t$ denotes the quotient algebra $\g/[\g,\g]_{\Q,L}$ and $\Lambda_T$ the image of the lattice $\Lambda$ in $\t$, then the Albanese variety of $X$ is nothing but $\t/\Lambda_T$ (see Subsection~\ref{algebraic_dimension_section}). It follows that the Albanese variety of $X$ is positive dimensional.

\hfill

{\bf Proof of (ii):} Let $Z\subset G$ be the center of $G$. Define $G_1$ to be the quotient group $G/Z$ and $\Gamma_1$ to be the image of $G$ in $G_1$ by the quotient map. Consider the map
    \begin{equation}
        \pi\colon X = \Gamma\backslash G \to X_1:= \Gamma_1\backslash G_1
    \end{equation}
    By \ref{_principal_toric_fiber_bundle_lemma_} the map $\pi$ is a principal toric fiber bundle. We repeat this construction with $X_1 = \Gamma_1\backslash G_1$ and so on to obtain the desired sequence of maps.
\endproof{}

\hfill


There exists a similar statement for abelian hypercomplex nilmanifolds. We will omit its proof as it follows the same lines as the proof of \ref{albanese_non_trivial_cor}.

\hfill

\proposition
\label{albanese_non_trivial_hypercomplex_prop}
Let $X = \Gamma\backslash G$ be an abelian hypercomplex nilmanifold (\ref{hypercomplex_nilmanifold_def}). Then
\begin{description}
    \item[(i)] The $\H$-Albanese variety (\ref{_h_albanese_def_}) of $X$ is positive-dimensional.
    \item[(ii)] There exists a finite sequence of submersions of hypercomplex manifolds
    \begin{equation}
        X = X_0 \to X_1 \to...\to X_{N-1}\to X_N = \{pt\}
    \end{equation}
    This sequence satisfies the following properties. First, for every $i$ the manifold $X_i$ is an abelian hypercomplex nilmanifold and the maps $X_i\to X_{i+1}$ are morphisms of hypercomplex manifolds. Second, the map $(X_i)_L\to (X_{i+1})_L$ is a holomorphic principal toric fiber bundle for any complex structure $L\in\H$.
\end{description}
\endproof{}


\section{Subvarieties in abelian hypercomplex nilmanifolds}

\subsection{Subvarieties in hypercomplex tori and nilmanifolds}

Let $X = \Gamma\backslash G$ be a nilmanifold and $\g$ the Lie algebra of $G$. Consider the smooth trivialization of the tangent bundle to $G$ by the left action of $G$ on itself. For every $g\in G$ we identify $\g$ with $T_gG$ by means of the operator $L_g$ of the left multiplication on $g$:

\begin{equation}
\label{_trivialization_eq_}
    (L_g)_*\colon \g \to T_gG
\end{equation}

This trivialization descends to a smooth trivialization of the tangent bundle to $X = \Gamma\backslash G$. Indeed, the left action of $\Gamma$ on $G$ preserves the trivialization by left-invariant vector fields. Hence we may and will identify $T_xX$ and $\g$ for every point $x\in X$.

\hfill

\definition
\label{_locally_homogeneous_submanifold_def_}
Let $X = \Gamma\backslash G$ be a nilmanifold and $Y\subset X$ be a submanifold of $X$. The submanifold $Y$ is said to be a {\em locally homogeneous submanifold of $X$} if there exists a subalgebra $\h\subset \g$ such that for every point $y\in Y$ the tangent space $T_yY$ is identified with $\h\subset\g$ via (\ref{_trivialization_eq_}).

\hfill 

Let $G$ be a Lie group with a left-invariant complex structure. By construction, the complex structure on $G$ is preserved by the map (\ref{_trivialization_eq_}). However, the trivialization of $TG$ via (\ref{_trivialization_eq_}) is not holomorphic unless $G$ is a complex Lie group. Let us explain this phenomenon in more detail. For any vector $\xi\in\g$ let $\tilde{\xi}$ denote the corresponding left-invariant vector field on $G$. A vector field on a smooth manifold is holomorphic if and only if its flow acts by holomorphic automorphisms. The flow of a {\em left-invariant} vector field $\tilde{\xi}$ is the multiplication by $\exp(t\xi)$ {\em on the right}. Hence the trivialization of $TG$ is holomorphic if and only if the complex structure is bi-invariant.

If $X=\Gamma\backslash G$ is a complex nilmanifold then the identification of $T_xX$ with $\g$ preserves the complex structure. The constructed trivialization of $TX$ is generally not holomorphic. It is holomorphic if and only if $G$ is a complex Lie group. If $X = \Gamma\backslash G$ is a hypercomplex nilmanifold then the identification of $T_xX$ with $\g$ preserves the hypercomplex structure on these spaces.

The main goal of this section is to give a proof of the following theorem.

\hfill

\theorem
\label{my_precious}
Let $X = \Gamma\backslash G$ be an abelian hypercomplex
nilmanifold. Then for all but at most countable number of
complex structures $L\in\H$, the complex nilmanifold $X_L$
satisfies the following property: if $Z\subset X_L$ is an
irreducible complex subvariety of $X_L$, then $Z$ is a
trianalytic locally homogeneous submanifold of $X$.

\hfill

The proof of this theorem will be given in the next subsection. Now we state its corollary, which may be seen as an analogue of \ref{submanifolds_in_hyperkahler_prop} for nilmanifolds with an HKT structure (\ref{HKT_def}).

\hfill

\corollary
\label{submanifolds_in_HKT_nilmanifold_trianalytic_cor}
Let $X$ be a hypercomplex nilmanifold with an HKT structure. Then for all but countable number of complex structures $L\in\H$, the complex nilmanifold $X_L$ satisfies the following property: if $Z\subset X_L$ is an irreducible complex subvariety of $X_L$, then $Z$ is a trianalytic submanifold of $X$ (\ref{trianalytic_def}).

\hfill

\proof By \ref{hypecomplex_abelian_equals_HKT} a hypercomplex nilmanifold $X$ admits an HKT metric if and only if the hypercomplex structure on $X$ is abelian. The result now follows from \ref{my_precious}. 

\hfill

\ref{my_precious} is well known in the case when $X$ is a hypercomplex torus.

\hfill

\proposition (\cite[Lemma 6.3]{Kaledin_Verbitsky})
\label{subvarieties_in_tori_prop}
Let $T = V/\Lambda$ be a hypercomplex torus, i.e. a quotient of a quaternionic vector space $V$ by a maximal rank lattice $\Lambda$. For all but at most countable number of complex structures $L\in\H$ the complex torus $T_L$ satisfies the following property: if $Z\subset T_L$ is an irreducible complex subvariety of $T_L$ then $Z$ is a hypercomplex subtorus.

\hfill

\proof We endow $T$ with a flat hyperk\"ahler metric. \ref{submanifolds_in_hyperkahler_prop} implies that for all but a countable number of complex structures $L\in \H$ every complex subvariety $Z\subset T_L$ is trianalytic. A trianalytic subvariety $Z$ of a hyperk\"ahler manifold is totally geodesic by \cite[Cor. 5.4]{Verbitsky_deformations_of_trianalytic}. The preimage $Z'$ of $Z$ in the vector space $V$ has to be an affine subspace. Indeed, a totally geodesic subvariety of a vector space is an affine subspace. Hence $Z$ is a subtorus.

\hfill

\remark Though some results of \cite{Kaledin_Verbitsky} are wrong (\cite{_Kaledin_Verbitsky:erratum_}), \cite[Lemma 6.3]{Kaledin_Verbitsky} is correct.

\hfill






\subsection{Subvarieties in principal toric fiber bundles}

Let $X = \Gamma\backslash G$ be a complex nilmanifold with a complex structure $L$, $\g$ the Lie algebra of $G$ and $\z$ (resp. $Z$) the center of the Lie algebra $\g$ (resp. the center of the Lie group $G$). Consider a holomorphic map
$$
\pi\colon \Gamma\backslash G=X \arrow Y:=\frac{G/Z}{\Gamma/(\Gamma\cap Z)}
$$
This map is principal toric bundle by \ref{_principal_toric_fiber_bundle_lemma_}. 

\hfill

\lemma
\label{semicontinuity_theorem}
Let $M\subset X$ be an irreducible complex subvariety,
for any $y\in Y$ denote $M_y:= \pi^{-1}(y)\cup M$.
Assume that $M_y\subset \pi^{-1}(y)$ is a finite union
of subtori in a compact torus $\pi^{-1}(y)$.
Given $x\in M_y$ we denote by $\z_x$ the Zariski
tangent space to $M_y$ at $x$ in the torus $\pi^{-1}(y)$.
\begin{description}
    \item[(i)] The family of subspaces $(\mathfrak z_x), x\in M$ of $\g$ is constant on a dense open subset $M'$ of $M$.
    \item[(ii)] Define $\mathfrak z':=\mathfrak z_x$ for $x\in M'$ and $Z':=\exp(\z')$. Then the subvariety $M$ is preserved by the right action of $Z'$.
\end{description} 

\hfill

{\bf Proof of (i):} Let $Y'$ be the maximal open subset of $Y$ such that $M':=\pi^{-1}(Y')\cap M$ contains only smooth points of $M$ and the map $\pi\restrict{M'}\colon M'\to Y'$ is smooth. The dimension of $\z_x$ is constant for $x\in M'$ and the subspace $\z_x$ depends continuously on $x$. Since the subspace $\z_x$ is rational for every $x\in M$, the family $\z_x$ is locally constant on $M'$. The subset $Y'$ is the complement to a closed analytic subset of $Y$. Hence $M'=\pi^{-1}(Y')\cap M$ is the complement to a closed analytic subset of $M$. Since $M$ is irreducible, $M'$ is connected. Therefore, the family $\z_x$ is constant on $M'$.

\hfill

{\bf Proof of (ii):} The right action of $Z'=\exp(\z')$ on $X$ preserves $M'$. Indeed, for any $y\in Y'$ the subvariety $\pi^{-1}(y)\cap M$ is a finite union of orbits of $Z'$. Since $M'$ is dense in $M$, the subvariety $M$ must be preserved by the action of $Z'$.
\endproof{}

\hfill

\subsection{Preliminary lemmas from algebraic geometry and topology}

We recall two well known results of algebraic geometry to be used in the sequel.

\hfill

\proposition
\label{fibers_are_divisors_prop}
Let $f\colon X\to Y$ be a birational morphism of complex varieties. Assume that $Y$ is smooth. Define the exceptional set $E\subset X$ as the union of positive dimensional fibers of $f$. Then $E$ is a divisor in $X$.

\hfill

\proof The statement is local on $Y$, hence we may assume that there exists a nowhere vanishing section $\omega$ of the canonical bundle of $Y$. The form $f^*\omega$ is a section of the canonical bundle of $X$. The zero set of $f^*\omega$ is precisely the exceptional set $E$. Thus $E$ is a (Cartier) divisor. \endproof{}

\hfill

The proof of the following proposition is almost the same.

\hfill

\proposition 
\label{ramification_locus_is_a_divisor}
Let $f\colon X\to Y$ be a finite morphism of complex varieties. Assume that $Y$ is smooth. Then the branch locus of $f$ is a divisor in $Y$.

\hfill

\proof The statement is local on $Y$, hence we may assume that there exists a nowhere vanishing section $\omega$ of the canonical bundle of $Y$. The form $f^*\omega$ is a section of the canonical bundle of $X$. The zero set of $f^*\omega$ is precisely the ramification locus $R\subset X$. Thus $R$ is a (Cartier) divisor in $X$. The branch locus is the image of $R$ in $Y$, hence a divisor in $Y$. \endproof{}

\hfill

We finish this subsection with a lemma of a topological nature. It will be used in the last step of the proof of \ref{my_precious}.

\hfill

\definition Let $\pi:\; X \arrow Y$ be a submersion of complex varieties.
Recall that {\bf a multisection} of $\pi$ is a subvariety
$Y_1 \subset Y$ such that the restriction $\pi:\; Y_1 \arrow Y$
is a finite map.

\hfill

\lemma
\label{multisections_lemma}
Let $X=\Gamma\backslash G$ be a nilmanifold,
$Z\subset G$ be the center of $G$. Consider the map
\begin{equation}
    \sigma\colon X = \Gamma\backslash G \to \frac{ G/Z} {\Gamma/(\Gamma\cap Z)}=: Y
\end{equation}
Then $\sigma$ admits no multisections.

\proof We use $\sigma_*$ to denote the map $\pi_1(X)\to \pi_1(Y)$ induced by $\sigma$. Notice that the center of $\pi_1(X)$
is precisely $\ker \sigma_*$ because
$Y= \frac{G/Z} {\Gamma/(\Gamma\cap Z)}$, and $Z$ is the center of $G$.

First, notice that
$\sigma$ admits no sections.
Indeed, $\sigma$ is a principal bundle,
hence any section of $\sigma$ induces a smooth decomposition
$X = T^n\times Y$, with $T^n = Z/(\Gamma\cap Z)$ a torus; here, $\sigma$ being the projection to the
second component. Then $\pi_1(X) = \pi_1 (T^n)\times \pi_1(Y)$,
with $\pi_1(T^n)$ being the center of $\pi_1(X)$. This is impossible;
indeed, since $Y$ is a nilmanifold, the center of $\pi_1(Y)$ has
positive rank, hence the center of $\pi_1(X) = \pi_1 (T^n)\times \pi_1(Y)$
is strictly bigger than $\pi_1(T^n)$: a contradiction.

To pass from multisections to sections,
consider the following operation on principal
toric bundles. Let $\pi\colon P \arrow V$ be a principal
toric bundle with fiber $T^n$, and
$P^k:= \underbrace{P\times_\sigma P \times_\sigma ... \times_\sigma
P}_{\text{$k$ times}}$
be its $k$-th power (as a principal bundle).
Taking a quotient of $P^k$ by a relation
$(p_1+t_1, p_2+t_2, ..., p_i+t_i, ..., p_k+t_k)\sim (p_1, p_2, ..., p_i,
..., p_k)$,
for all $t_1,t_2,...t_k\in T^n$ such that $\sum t_i = 0$, we obtain a principal
bundle $P_k$ with a fiber $T^n$ again. However, this bundle
is not, generally speaking isomorphic to $P$; the natural map
$P \arrow P_k$ mapping $p$ to $(p, p, ..., p)$ is a $k$-torsion
quotient of $P_k$. We call the bundle $P_k$ {\bf the
quotient of $P$ by $k$-torsion}. Any order $k$ multisection
of $P$ gives a section of $P_k$, hence trivializes $P_k$.

Return now to the proof of \ref{multisections_lemma}.
A multisection of $\sigma$ gives a section of
the principal bundle $\sigma_k\colon X_k \arrow Y$,
where $X_k$ is the quotient of the principal
bundle $\sigma\colon X \arrow Y$ by $k$-th torsion.
Indeed, let $\rho$ be an order $k$ multisection mapping
$y\in Y$ to $(s_1, ..., s_k)$. Then
$(s_1, ..., s_k)$ defines a point in $X_k$ which
is independent from the ordering of $(s_1, ..., s_k)$.

However, $X_k$ is also a nilmanifold, with the
same rational Lie algebra (indeed, a nilmanifold
can be defined as a compact manifold which is
homogeneous under an action of a nilpotent
Lie group).  Therefore, the projection $X_k\arrow Y$
corresponds to taking the quotient of the corresponding
Lie group by its center, and 
$\ker \left((\sigma_k)_*\colon \pi_1(X_k)\to \pi_1(Y)\right)$
is the center of $\pi_1(X_k)$. Since the center of a
nilpotent Lie algebra cannot be its direct summand,
the projection $\sigma_k$ also cannot have sections. 
We proved \ref{multisections_lemma}.
\endproof


\subsection{Subvarieties of hypercomplex nilmanifolds: the proof}
\label{_Proof_Subsection_}.

This subsection is completely dedicated to the proof of \ref{my_precious}. 

\hfill

{\bf Theorem 5.2:} Let $X = \Gamma\backslash G$ be an
abelian hypercomplex nilmanifold. Then for all but at most
countable number of complex structures $L\in\H$, the
complex nilmanifold $X_L$ satisfies the following
property: if $Z\subset X_L$ is an irreducible complex
subvariety of $X_L$, then $Z$ is a trianalytic locally
homogeneous submanifold of $X$.

\hfill

The proof goes by induction on $n = \dim_\H X$. It is divided into steps for the convenience of the reader. In {\bf Step 0} we prove the base of induction. In {\bf Step 1} we consider the principal toric fiber bundle $\pi\colon X\to Y$ as in the second part of \ref{albanese_non_trivial_hypercomplex_prop}. We remind the construction of this map and prove the assertion of the theorem in the case then $Y$ is a point. In {\bf Step 2} we define the set of ``nice" complex structures on $X$ and prove that this set is a complement to at most countable number of complex structures in $\H$. In the next steps we will prove that for every ``nice" complex structure $L$ each irreducible complex subvariety of $X_L$ is a hypercomplex locally homogeneous submanifold of $X$. In {\bf Step 3} we show that it is enough to prove the theorem for complex subvarieties of $X_L$ which are mapped surjectively onto $Y$. In {\bf Step 4} we reduce the theorem further: it is enough to check the statement of the theorem for complex subvarieties $M\subset X_L$ such that the map $\pi\restrict{M}\colon M\to Y$ is surjective and generically finite. In {\bf Step 5} we show that if the map $\pi\restrict{M}\colon X\to Y$ is generically finite then it is \'etale. It is impossible for $\pi\restrict{M}$ to be \'etale by \ref{multisections_lemma}. This finishes the proof.

\hfill

{\bf Proof of \ref{my_precious}. Step 0:} The proof goes by induction on $n = \dim_\H X$. The only abelian hypercomplex nilmanifold $X =\Gamma\backslash G$ of real dimension $4$ is a torus (see also \cite{_Boyer_},\cite{_Hasegawa_}). Indeed, the center of the Lie algebra $\g$ must be $\H$-invariant by \ref{upper_central_series_L_invariant_lemma}, hence coincide with $\g$. The statement of the theorem for a hypercomplex torus of real dimension $4$ follows from \ref{subvarieties_in_tori_prop} or directly from \ref{submanifolds_in_hyperkahler_prop}.

\hfill

{\bf Step 1:} Assume now that $\dim_\H X = n>1$ and the statement of the theorem holds for all abelian hypercomplex nilmanifolds of quaternionic dimension less than $n$. Consider the morphism ${\pi\colon X\to Y}$ of hypercomplex manifolds as in the second part of \ref{albanese_non_trivial_hypercomplex_prop}. We briefly recall its construction here. Let $\z$ (resp. $Z$) denote the center of the Lie algebra $\g$ (resp. the center of the Lie group $G$). Then the quotient map ${G\to G/Z}$ induces the map ${\pi\colon \Gamma\backslash G = X \to Y:= \frac{G/Z}{\Gamma/(\Gamma\cap Z)}}$. For every complex structure $L\in\H$ the map $\pi$ is a principal holomorphic toric fiber bundle whose fibers are isomorphic to the hypercomplex torus ${T:= Z/(\Gamma\cap Z)}$ (\ref{albanese_non_trivial_hypercomplex_prop}).

In the case when the manifold $Y$ is a point, the manifold $X$ is a torus. The assertion in this case follows from \ref{subvarieties_in_tori_prop}.

\hfill

{\bf Step 2:} By Step 1 we may assume that the manifold ${Y = \frac{G/Z}{\Gamma/(\Gamma\cap Z)}}$ has positive dimension. We introduce the following property of a hypercomplex nilmanifold $N$ and a complex structure $L\in\H$.
\begin{equation}
\label{_property_}
    \begin{minipage}[t]{0.65\linewidth}
    {\em Every irreducible complex subvariety of $N_L$ is a trianalytic locally homogeneous submanifold.}
    \end{minipage}
\end{equation}

Let us call a complex structure $L\in \H$ ``nice" if it satisfies the following two properties.

\begin{enumerate}

\item Take a proper $\H$-invariant rational subalgebra $\g'$ of the Lie algebra $\g$. Let $G'$ denote the exponent of $\g'$ in $G$. Clearly, the quotient space $X':=(\Gamma\cap G')\backslash G'$ is a hypercomplex nilmanifold embedded into $X$ as a trianalytic locally homogeneous submanifold. Then for every choice of $\g'$ the nilmanifold $X'$ satisfies \eqref{_property_}.

\item Take an $\H$-invariant rational subspace $\z'$ of the center $\z$ of $\g$. Let $Z'$ denote the exponent of $\z'$ in $G$. Clearly, the quotient space ${Y':= \frac{G/Z'}{\Gamma/(\Gamma\cap Z')}}$ is an abelian hypercomplex nilmanifold. Then for every choice of $\z'$ the nilmanifold $X'$ satisfies \eqref{_property_}.

\end{enumerate}

We claim that the set of ``nice'' complex structures $L\in \H$ is the complement to at most countable number of complex structures. Indeed, the set of rational subspaces of $\g$ is countable. It follows from the induction hypothesis that the first and the second condition above each exclude at most countable number of complex structures.

In the next steps we will prove that every complex subvariety of $X_L$ is trianalytic as soon as the complex structure $L$ is ``nice".

\hfill

{\bf Step 3:} Take a nice complex structure $L\in\H$ and let $M\subset X_L$ be an irreducible complex subvariety of complex dimension $k$. Consider the image $\pi(M)$ of $M$ under the map ${\pi\colon X\to Y = \frac{G/Z}{\Gamma/(\Gamma\cap Z)}}$ (see Step 1). Suppose that $\pi(M)$ is properly contained in $Y$. The second condition on a ``nice" complex structure applied to $\z'=\z$ implies that $\pi(M)$ is a trianalytic locally homogeneous submanifold of $Y$. By construction, the preimage $X':=\pi^{-1}(\pi(M)))$ of $\pi(M)$ in $X$ is a trianalytic locally homogeneous submanifold of $X$. It is associated with some rational $\H$-invariant subalgebra $\g'\subset\g$. By the first condition on a ``nice'' complex structure every complex subvariety of $X'_L$ is a trianalytic locally homogeneous submanifold. Since $M$ is contained in $X'$, it is trianalytic locally homogeneous. We have reduced the theorem to the case when $M$ is mapped surjectively onto $Y$.

\hfill

{\bf Step 4:} By Step 3 we may assume that $M$ is mapped surjectively onto $Y$. Pick a point $y\in Y$. Consider the fiber $\pi^{-1}(y)\cap M$ of the map $\pi\restrict{M}\colon M\to Y$. The first condition on a ``nice'' complex structure applied to $\g' = \z$ implies that $\pi^{-1}(y)\cap M$ is a trianalytic locally homogeneous submanifold of $\pi^{-1}(y)$ and hence of $X$. For every $x\in\pi^{-1}(y)$ the tangent space $T_x(\pi^{-1}(y)\cap M)$ is identified with a rational $\H$-invariant subspace $\z_x\subset \z$ via the trivialization (\ref{_trivialization_eq_}). By \ref{semicontinuity_theorem} there exists an $\H$-invariant rational subspace $\z'\subset\z$ such that $\z'\subset \z_x$ for every $x\in M$ and $\z'=\z_x$ for $x$ in a dense open subset of $X$. Moreover, the right action of the group $Z':=\exp(\z')$ preserves $M$.

Suppose that $\dim \z'>0$. Consider the map 
\begin{equation}
    \pi'\colon X \arrow Y':=\frac{G/Z'}{\Gamma/(\Gamma\cap Z')}
\end{equation}
The map $\pi'$ is a holomorphic principal toric fiber bundle by \ref{_principal_toric_fiber_bundle_lemma_}. By the second condition on a ``nice'' complex structure, $\pi'(M)$ is a trianalytic locally homogeneous submanifold of $Y'$. Since $M = (\pi')^{-1}(\pi'(M))$, it is a trianalytic locally homogeneous submanifold of $X$. We reduced the theorem to the case when $\dim \z' = 0$, or in other words, the map $\pi\restrict{M}\colon M\to Y$ is surjective and generically finite.

\hfill

{\bf Step 5:} Consider a subvariety $M\subset X$ such that the map $\pi\restrict{M}\colon M\to Y$ is generically finite. Consider the Stein factorization of the map $\pi\restrict{M}$:
$$
\begin{CD}
    M @>{\pi_0}>> Y'@>{\pi_1}>> Y
\end{CD}
$$
Here $\pi_0\colon M\to Y'$ has connected fibers and $\pi_1\colon Y'\to Y$ is a finite morphism. By \ref{ramification_locus_is_a_divisor} the branch locus of $\pi_1$ is a divisor in $Y$. Every complex subvariety of $Y_L$ is trianalytic, in particular, there are no divisors in $Y_L$. Hence the morphism $\pi_1\colon Y'\to Y$ is \'etale. In particular, $Y'$ is smooth. The map $\pi_0\colon M\to Y'$ is a birational map to a smooth base. By \ref{fibers_are_divisors_prop} the exceptional set $E\subset M$ is a divisor in $M$. If $E$ is non-empty then $E$ has odd dimension. Consider the image $\pi(E)$ of $E$ in $Y$. It is a complex subvariety of $Y_L$, hence trianalytic. In particular, $\pi(E)$ is of even dimension. Therefore, for a generic point $y\in \pi(E)$ the variety $\pi^{-1}(y)\cap E$ is odd-dimensional. All subvarieties of $\pi^{-1}(y)$ are trianalytic and hence of even dimension. Hence $E$ must be empty. We have shown that $\pi_0\colon M\to Y'$ is an isomorphism and $\pi_1\colon Y'\to Y$ is \'etale. Hence the morphism $\pi\restrict{M}\colon M
\to Y$ is \'etale. By \ref{multisections_lemma} there are no subvarieties $M\subset X$ with this property. The theorem is proved.
\endproof{}




\hfill


\hfill

{\bf Acknowledgements:}
We are grateful to Nikon Kurnosov for pointing out
to us the results of \cite{_Bogomolov_Kurnosov_Kuznetsova_Yasinsky_} and fruitful
discussions, Dmitry Krekov and Andrey Soldatenkov for their invaluable help and
Iuliia Gorginian for her insight and comments.

\hfill

{\small 
\noindent {\sc Anna Abasheva\\
{\sc Columbia University\\
Department of Mathematics, 2990 Broadway, New York, NY, USA}\\
also:\\
\sc Independent University of Moscow\\
Bolshoy Vlasievskiy per., 11, Moscow, Russia\\
		also:\\
		Laboratory of Algebraic Geometry, \\
		Faculty of Mathematics, HSE University,\\
		7 Vavilova Str. Moscow, Russia}\\
\tt  anabasheva@yandex.ru

\hfill

	\noindent \sc Misha Verbitsky\\
		\sc Instituto Nacional de Matem\'atica Pura e
			Aplicada (IMPA) \\ Estrada Dona Castorina, 110\\
			Jardim Bot\^anico, CEP 22460-320\\
			Rio de Janeiro, RJ - Brasil \\
		also:\\
		Laboratory of Algebraic Geometry, \\
		Faculty of Mathematics, HSE University,\\
		7 Vavilova Str. Moscow, Russia\\
	\tt verbit@impa.br }
\end{document}